\documentclass{article}
\usepackage{latexsym,amsfonts,amsmath,amsthm,amssymb,makeidx}
\usepackage[title]{appendix}
\usepackage{CJK,CJKnumb,CJKulem,times,dsfont,ifthen,mathrsfs,latexsym,amsfonts, color}
\usepackage{amsmath,amsthm,makeidx,fontenc,amssymb,bm,graphicx,psfrag,listings, curves,extarrows}

\let\oldbibliography\thebibliography
\renewcommand{\thebibliography}[1]{%
\oldbibliography{#1}%
\setlength{\itemsep}{0pt}%
}

\makeindex
\newtheorem{definition}{Definition}[section]
\newtheorem{theorem}{Theorem}[section]
\newtheorem{lemma}{Lemma}[section]
\newtheorem{corollary}{Corollary}[section]

\newtheorem{remark}{Remark}[section]

\newcommand{\q}{\theta}

\newcommand{\la}{\lambda}

\newcommand{\pa}{\partial}

\newcommand{\R}{\mathbb R}

\newcommand{\C}{\mathbb C}

\newcommand{\bt}{\begin{theorem}}
\newcommand{\et}{\end{theorem}}
\newcommand{\bl}{\begin{lemma}}
\newcommand{\el}{\end{lemma}}
\newcommand{\bd}{\begin{definition}}
\newcommand{\ed}{\end{definition}}
\newcommand{\bc}{\begin{corollary}}
\newcommand{\ec}{\end{corollary}}
\newcommand{\bp}{\begin{proof}}
\newcommand{\ep}{\end{proof}}
\newcommand{\bx}{\begin{example}}
\newcommand{\ex}{\end{example}}
\newcommand{\bi}{\begin{exercise}}
\newcommand{\ei}{\end{exercise}}
\newcommand{\bo}{\begin{prop}}
\newcommand{\eo}{\end{prop}}
\newcommand{\br}{\begin{remark}}
\newcommand{\er}{\end{remark}}
\newcommand{\be}{\begin{equation}}
\newcommand{\ee}{\end{equation}}
\newcommand{\ba}{\begin{align}}
\newcommand{\ea}{\end{align}}
\newcommand{\bn}{\begin{enumerate}}
\newcommand{\en}{\end{enumerate}}
\newcommand{\bg}{\begin{align*}}
\newcommand{\bcs}{\begin{cases}}
\newcommand{\ecs}{\end{cases}}

\newcommand{\md}{{\mathbf {div}_{\theta}}}
\newcommand{\ds}{{\nabla_{\theta}}}

\newcommand{\bean}{\begin{eqnarray*}}
\newcommand{\eean}{\end{eqnarray*}}

\numberwithin{equation}{section}

\begin{document}

\title{\bf On finite Morse index solutions to the  quadharmonic  Lane-Emden equation \thanks{Partially supported by NSFC of China and  NSERC of Canada. E-mails: luosp14@mails.tsinghua.edu.cn(Luo);\;  jcwei@math.ubc.ca (Wei);\; wzou@math.tsinghua.edu.cn(Zou)} }
\date{}
\author{\\{\bf  Senping Luo$^{1}$,\;\;  Juncheng Wei$^{2}$ \;\; and \; Wenming Zou$^{3}$}\\
\footnotesize {\it  $^{1,3}$Department of Mathematical Sciences, Tsinghua University, Beijing 100084, China}\\
\footnotesize {\it  $^{2}$Department of Mathematics, University of British Columbia, Vancouver, BC V6T 1Z2,
Canada}
}

\maketitle
\begin{center}
\begin{minipage}{120mm}
\begin{center}{\bf Abstract}\end{center}

In this paper, we compute the Joseph-Lundgren exponent for the  quadharmonic Lane-Emden equation, derive a   monotonicity formula and classify  the finite Morse index solution  to the following quadharmonic Lane-Emden equation:
\noindent
\begin{equation}\nonumber
\Delta^4 u=|u|^{p-1}u\;\;\;\;\hbox{in}\;\;\;\;\; \R^n.
\end{equation}
As a byproduct, we also get a monotonicity formula for the quadharmonic maps $ \Delta^4 u=0$.

\vskip0.10in


\end{minipage}
\end{center}

\vskip0.10in
\section{Main results and Background}

 We study the finite Morse index solution of the following Lane-Emden equation
 \begin{equation}\label{12LE}
\Delta^4 u=|u|^{p-1}u\;\;\;\;\hbox{in}\;\;\;\;\;  \R^n.
\end{equation}
The goal of the current paper is to complete the classification of the finite Morse index solution to Eq. \eqref{12LE} and to establish the corresponding Liouville-type theorems. It is well known that the Liouville-type theorems  play a crucial role to obtain  a priori $L^\infty$-bounds for solutions of semilinear   elliptic and parabolic problems.  See the monograph of Quittner and Souplet \cite{QS2007}. To the best of our knowledge,
such a kind of work  to Eq. \eqref{12LE} has not been accomplished  previously.  Postponing  the background of this topic to the last part of this section, we would like to state the main results of the current article  first.  Recall  that a  solution $u$ of \eqref{12LE} is said to be    stable outside a compact $\Theta\subset \R^n $ if
\be\nonumber
\int_{\R^n}|\Delta^2 \varphi|^2dx\geq p\int_{\R^n}|u|^{p-1}\varphi^2dx,\quad\hbox{for any}\;\varphi\in H^4(\R^n\backslash\Theta).
\ee
In particular, if $\Theta=\emptyset$, we say that  $u$  is stable on $\R^n$. On the other hand, the Morse index   of the solution $u$ of \eqref{12LE} is defined as the maximal dimension over  all subspaces $E$ of $H^4 (\R^n)$ satisfying
\be\nonumber
\int_{\R^n}|\Delta^2 \varphi|^2dx<p\int_{\R^n}|u|^{p-1}\varphi^2dx,\quad\hbox{for any}\; \varphi\in E\setminus\{0\}.
\ee
 It is known that if a solution $u$ to \eqref{12LE} has finite Morse index,  then $u$ must be  stable outside a compact of $  \R^n. $
 The first main result  of the present paper is the following

\bt\label{12th1}
 Let $u$ be a stable solution of \eqref{12LE}.  If $ 1<p< p_c (n)$, then $ u \equiv 0$.
\et
Where $p_c(n):=p_{cQuadharmonic}(n)$ is the  Joseph-Lundgren exponent for the quadharmonic Lane-Emden equation  \eqref{12LE}, which
will be introduced shortly.   Further, we have more general results for  the finite Morse index solutions.

\bt\label{12th2}
Let $u$ be a finite Morse index solution of \eqref{12LE}.
Assume that either
\begin{itemize}
\item [(1)]   $1<p<\frac{n+8}{n-8}$
or
\item [(2)]  $\frac{n+8}{n-8}< p<p_c(n)$,
\end{itemize}
then the solution $u\equiv0$.
\begin{itemize}
\item[(3)]  If $p=\frac{n+8}{n-8}$, then $u$ has a finite energy, i.e.,
\be\nonumber
\int_{\R^n}|\Delta^2 u|^2=\int_{\R^n}|u|^{p+1}<+\infty.
\ee

\end{itemize}
\et

\br
In the above both Theorems,  $p_c (n)$, given below, called the Joseph-Lundgren exponent.  The condition $p <p_c (n)$ is optimal. In fact the radial singular solution is stable when $p\geq p_c (n)$ (See \cite{LWZ}).
\er
Now we give the Joseph-Lundgren exponent for the quadharmonic Lane-Emden equation,  which is defined by
\be\label{12pcn}
p_c(n):=p_{cQuadharmonic}(n):=
\begin{cases}
\;\;\;\infty\;\;&\hbox{if}\;\; n\leq 17,\\
\frac{n+6-2d(n)}{n-10-2d(n)}\;\;&\hbox{if}\;\; n\geq18,
\end{cases}
\ee
where
\be\label{LWZ=333}\aligned
d(n):=\sqrt{\frac{1}{4}n^2+5+\frac{1}{2}\sqrt{d_6}-\frac{1}{2}\sqrt{d_7+\frac{d_3}{\sqrt{d_6}}}}
\endaligned\ee
and
\be\nonumber\aligned
d_0:&=2097152-\frac{45}{4}n^{10}+180n^9-396n^8-5184n^7+36928n^6+27648n^5\quad\quad\quad\quad\quad \\
&\quad -132096n^4+147456n^3-1572864n^2;
\endaligned\ee
\be\nonumber\aligned
d_1:&=\frac{3}{65536}n^{24}-\frac{9}{4096}n^{23}+\frac{81}{2048}n^{22}-\frac{33}{128}n^{21}-\frac{123}{128}n^{20}
+\frac{303}{16}n^{19}+\frac{21}{8}n^{18}\\
&\quad-1056n^{17}+3888n^{16}+25396n^{15}-279456n^{14}+947712n^{13}+1979904n^{12}\\
&\quad-48427008n^{11}+135979008n^{10}+677117952n^9
-2620588032n^8\\
&\quad-3265265664n^7+14294188032n^6
+2415919104n^5-16106127360n^4;
\endaligned\ee
\be\nonumber\aligned
d_2:=(d_0+12\sqrt{d_1})^{\frac{1}{3}};d_3:=128n^2; \quad\quad\quad\quad\quad \quad\quad\quad\quad\quad \quad\quad\quad\quad\quad \quad\quad\quad\quad\quad
\endaligned\ee
\be\nonumber\aligned
d_4:=-\frac{8192}{3}+\frac{1}{32}n^8-\frac{1}{2}n^7+n^6+16n^5-\frac{584}{3}n^4-128n^3+\frac{4096}{3}n^2;\quad\quad\quad\quad\quad
\endaligned\ee
\be\nonumber\aligned
d_5:=\frac{40}{3}n^2+\frac{128}{3}, \quad
d_6:=\frac{1}{2}d_5+\frac{1}{6}d_2-\frac{d_4}{d_2},\quad
d_7:=d_5-\frac{1}{6}d_2+\frac{d_4}{d_2}. \quad\quad\quad\quad\quad
\endaligned\ee
The expression of $d(n)$ seems very complicated,
 however we have that $$d(n)<\sqrt{n} \;\; \hbox{  for } n\geq18; \quad \lim_{n\rightarrow\infty}\frac{d(n)}{\sqrt{n}}=1.$$

\br In the quadharmonic case,  the above Joseph-Lundgren exponent $p_c (n)$ actually satisfies a $8$-th order polynomial algebraic  equation. It is interesting that we can obtain  the explicit formulation  here.\er

\br Theorems  \ref{12th1}-\ref{12th2} are proved by Monotonicity Formula which we introduce in the next section.\er
\vskip0.223in


Now we describe the background and development for lane-Emden equations. It is well known that the   Lane-Emden equation
 \begin{equation}\label{LZ=0001}
 -\Delta  u=|u|^{p-1}u \;\;\;\;\hbox{in}\;\;\;\;\;  \R^n
\end{equation}
and its parabolic counterpart have played a key  role in the development of methods and applications of nonlinear PDEs in the last
decades.  The fundamental works on  Eq.\eqref{LZ=0001} are due to  \cite{Fowler=1,Fowler=2}.
Another ground-breaking  result on equation \eqref{LZ=0001} is the celebrated   Liouville-type theorem due to Gidas and Spruck \cite{Gidas-Sp=1981}, they assert  that the Eq. \eqref{LZ=0001} has no positive solution whenever $p\in (1, 2^\ast-1)$, where $2^\ast=2n/(n-2)$ if $n\geq 3$ and
 $2^\ast=\infty$ if $n\leq 2.$ However, if $p=2^\ast-1$  the Eq. \eqref{LZ=0001}  has a unique positive solution (up to translation  and rescaling)
 which is radial and explicit (see Caffarelli-Gidas-Spruck \cite{Caffarelli1989}). Since then  there has been  an extensive literature on such a  type of equations  or systems.   Among them,  the paper  \cite{Farina2007} by Farina  in 2007 (see also \cite{Farina2005}), the equation  \eqref{LZ=0001} is revisited  for $p>\frac{n+2}{n-2}$.  The author obtained  some  classification results and Liouville-type theorems  for $\C^2(\R^n)$ smooth solutions  including stable solutions, finite Morse index solutions, solutions which are stable
outside a compact set, radial solutions and non-negative solutions. The  results obtained in  \cite{Farina2007} were   applied to subcritical, critical and supercritical values of the exponent $p.$  Moreover,   the critical stability exponent $p_c(n)$ (Joseph-Lundgren exponent) is determined which   is larger than the classical critical exponent  $p_S=2^\ast-1$ in Sobolev imbedding theorems. Precisely, in  \cite{Farina2007},  the $p_c(n)$    is given by
\be
p_c(n):=p_{cHarmonic}(n):=\begin{cases}
\;\;\;\;\;\;\;\;\;\infty\;\;\;\;\;\;\;\;&\hbox{if }\;  n\leq10,\\
\frac{(n-2)^2-4n+8\sqrt{n-1}}{(n-2)(n-10)}\;\;\;\;\;\;\;\;&\hbox{if }\;  n\geq11
\end{cases}
\ee
which can be traced back to Joseph-Lundgren \cite{Joseph1972}.  It is proved that the $\C^2(\R^n)$ stable solution of Eq.\eqref{LZ=0001} is identically to zero if $p<p_c(n)$, while  Eq.\eqref{LZ=0001} admits a smooth positive, bounded, stable and radial
solution if $p\geq p_c(n)(n\geq 11)$. In some sense, the   Joseph-Lundgren exponent  $p_c(n)$ is a critical threshold for obtaing the  Liouville-type theorems for stable or finite Morse index solutions.  The proof of Farina involves a delicate use of Nash-Moser's iteration technique, which is a classical tool for regularity of second order elliptic operators and falls short for higher order operators.

 \vskip0.1in

 About the biharmonic equation
 \begin{equation}\label{12DS-zwm}
\Delta^2 u=|u|^{p-1}u\;\;\;\;\hbox{in}\;\;\;\;\;  \R^n,
\end{equation}
the corresponding Joseph-Lundgren exponent (Seee Gazzola and Grunau \cite{Gazzola2006}, 2006) is
\be
p_c:=p_{cBiharmonic}(n):=\begin{cases}
\;\;\;\;\;\;\;\;\;\;\infty\;\;\;\;\;\;\;\;&\hbox{if}\;  n\leq12,\\
\frac{n+2-\sqrt{n^2+4-n\sqrt{n^2-8n+32}}}{n-6-\sqrt{n^2+4-n\sqrt{n^2-8n+32}}}\;\;\;\;\;\;\;\;&\hbox{if}\;  n\geq13.
\end{cases}
\ee
Further, Davila, Dupaigne, Wang and Wei in \cite{Wei=2}
obtain  the Liouville-type theorem and give a complete characterization of all finite Morse index solutions (whether radial or not,whether positive or not). See  Lin \cite{Lin1998}  and Wei-Xu \cite{Wei1999}  where it is shown that
  any  nonnegative solution $u$  of Eq.\eqref{12DS-zwm} is $\equiv 0$
for $1 < p <\frac{n + 4}{n-2}.$

\vskip0.1in

Very recently, in our previous manuscript \cite{LWZ2016=3}, we derive a new  monotonicity formula and classify completely all the  finite Morse index solutions (positive or sign-changing, radial or not) to  the   triharmonic Lane-Emden equation:
\begin{equation}\label{12DS=LWZ=100}
(-\Delta)^3 u=|u|^{p-1}u\;\;\;\;\hbox{in}\;\;\;\;\;  \R^n,
\end{equation}
where  the corresponding Joseph-Lundgren exponent is determined  by    the following formula:
\be\nonumber
p_{cTriHarmonic}(n):=p_c:=
\begin{cases}
\;\;\;\infty\;\;&\hbox{if}\;\; n\leq 14,\\
\frac{n+4-2D(n)}{n-8-2D(n)}\;\;&\hbox{if}\;\; n\geq15,
\end{cases}
\ee
where
\be\nonumber
D(n):=\frac{1}{6}\Big(9n^2+96-\frac{1536+1152n^2}{d_0(n)}-\frac{3}{2}d_0(n)\Big)^{1/2};
\ee
\be\nonumber
D_0(n):=-(D_1(n)+36\sqrt{D_2(n)})^{1/3};
\ee
\be\nonumber\aligned
D_1(n):=-94976+20736n+103104n^2-10368n^3+1296n^5-3024n^4-108n^6;
\endaligned\ee
\be\nonumber\aligned
D_2(n):&=6131712-16644096n^2+6915840n^4-690432n^6-3039232n\\
&\quad+4818944n^3-1936384n^5+251136n^7-30864n^8-4320n^9\\
&\quad+1800n^{10}-216n^{11}+9n^{12}.
\endaligned\ee
Obviously, the exponent $p_c(n)$ becomes more and more complex along with the order's increasing.

\vskip0.1in
 On the other hand, we note  the nonlocal Lane-Emden equation:
 \begin{equation}\label{12DS}
(-\Delta)^s u=|u|^{p-1}u\;\;\;\;\hbox{in}\;\;\;\;\;  \R^n.
\end{equation}
When $0<s<1$ and $1<s<2$, the complete classification of finite Morese index solution to Eq. \eqref{12DS} has been finished  by
Davila, Dupaigne, Wei in \cite{Wei0=1} and Fazly, Wei in \cite{Wei1=2} respectively.

 \vskip0.1in

Unfortunately, so far ones do not have a generic approach to deal with the  general polyharmonic Lane-Emden  case
 \begin{equation}\label{12DS=LWZ=m}
(-\Delta)^m u=|u|^{p-1}u\;\;\;\;\hbox{in}\;\;\;\;\;  \R^n
\end{equation}
or polyharmonic map $
(-\Delta)^m u=0\; \hbox{in}\;  \R^n $  to  obtain  the Liouville-type theorem and  a complete characterization of all finite Morse index solutions.
This remains an interesting open problem.

 \vskip0.1in

Finally, we refer the readers to  J. Serrin, H. Zou \cite{Serrin-Zou-1998}, P. Souplet \cite{Souplet2009} and  E. Mitidieri \cite{Mitidieri-96} for the  Lane-Emden  systems and     X. F. Wang \cite{WangXF-98-TAMS}  for the corresponding
reaction-diffusion equations.

 \vskip0.13in

The paper is organized as follows: In Section 2, we introduce the monotonicity formula (Theorems 2.1 and 2.2). In Subsections 2.1-2.2, we give
some preliminary calculations related to the  functional of the  monotonicity formula. In Section 3, we  give the representations
on the operators $\Delta^j,j=1,2,3$. In Section 4, we establish the differential by part formulas. In Sections 5-6, we calculate the derivatives
of the functional of the  monotonicity formula and prove Theorem 2.1.  The Section 7 is devoted to prove the desired monotonicity formula, i.e.,  Theorem \ref{12mono2}. In Section 8, we will show that the homogeneous stable solution must be zero. The Section 9 is on the energy estimates and blow-down analysis, we will prove Theorem 1.1. Finally, the Section 10 will study  the  finite Morse index solution and prove Theorem 1.2.


\section{Monotonicity formula }

\vskip0.1in

Let $(r,\theta)$ be the spherical coordinates in $\R^n,$  i.e., $ r = |x| \in  (0,\infty)$  and
$\theta = x/|x| $ is a point of the unit sphere $ S^{n-1}$.  The symbols $\Delta_\theta$ and $\nabla_\theta$ refer  respectively  to the Laplace-Beltrami operator and the gradient on  $ S^{n-1}$.
We denote $\pa_r u= \nabla u\cdot\frac{x}{r}, r=|x|$. Let
$$B_\lambda:=\{y\in \R^n: |y-x|<\lambda\}  $$
and
$$
 u^\la(x):=\la^{\frac{8}{p-1}}u(\la x),\;\;  \lambda>0. $$  Define
\be\label{12Exu}\aligned
&E(\la,x,u):=
\int_{B_1}\frac{1}{2}|\Delta^2 u^\la|^2-\frac{1}{p+1}|u^\la|^{p+1}\\
&\;\;+\int_{\pa B_1}\Big(\sum_{i,j\geq0,i+j\leq7}C_{i,j}^{0}\la^{i+j}\frac{d^iu^\la}{d\la^i}\frac{d^ju^\la}{d\la^j}
+\sum_{i,j\geq0,i+j\leq5}C_{i,j}^{1}\la^{i+j}\ds\frac{d^i u^\la}{d\la^i}\ds\frac{d^ju^\la}{d\la^j}\\
&\;\;+\sum_{i,j\geq0,i+j\leq3}C_{i,j}^{2}\la^{i+j}\Delta_{\theta}\frac{d^iu^\la}{d\la^i}\Delta_{\q}\frac{d^ju^\la}{d\la^j}\\
&\;\;+\sum_{i,j\geq0,i+j\leq1}C_{i,j}^{3}\la^{i+j}\ds\Delta_{\q}\frac{d^iu^\la}{d\la^i}\ds\Delta_{\q}\frac{d^ju^\la}{d\la^j}\Big).\\
\endaligned\ee

\bt\label{12mono1} Suppose that $u$ is a solution of \eqref{12LE}, then we have the following monotonicity formula
\be\label{12monoi}\aligned
&\frac{d}{d\la}E(\la,x,u)=
\int_{\pa B_1}\sum_{j=1}^4(A_j+a_j)\la^{2j-1}(\frac{d^ju^\la}{d\la^j})^2+2\la\int_{\pa B_1}\big(\nabla_\theta\Delta_{\theta} \frac{du^\la}{d\la}\big)^2\\
&+\int_{\pa B_1}\sum_{s=1}^3(B_s+b_s)\la^{2s-1}(\nabla_\theta\frac{d^su^\la}{d\la^s})^2
+\int_{\pa B_1}\sum_{l=1}^2(C_l+c_l)\la^{2l-1}(\Delta_\theta\frac{d^l u^\la}{d\la^l})^2\\
&+2\int_{\pa B_1}\la|\nabla_{\q}
\frac{dv^\la}{d\la}|^2,
\endaligned\ee
where the constants $C^s_{i,j}$ depending on $n,p$ can be determined in our proofs below.

 \et

\br In our proof below, the constants $C^s_{i,j}$ (and the constants in the next Corollary) can be determined via $n,p$ and may be negative,
but we don't need the exact expressions of them. Furthermore, we don't need the positiveness of the constants in our proof of the Liouville type theorem. Therefore,  we just represent them by  the generic notations  $C^s_{i,j}$.
\er

\bt\label{12mono2}
Suppose that $u$ is a solution of \eqref{12LE}. If $\frac{n+8}{n-8}<p<p_c(n)$, then there exists a constant $C(n,p)>0$
such that
\be\label{12monoi}\aligned
\frac{d}{d\la}E(\la,x,u)
&\geq C(n,p)\int_{\pa B_1}\la(\frac{du^\la}{d\la})^2\\
&=C(n,p)\la^{8\frac{p+1}{p-1}-8-n}\int_{\pa B_\la}(\frac{8}{p-1}+\la\pa_r u)^2
\endaligned\ee

\et

  The proof of Theorem \ref{12mono1} will be  split   in the Sections 2-6.
Based on  Theorem \ref{12mono1}, by combining the  algebraic and differential analysis in the Section 7, we can get Theorem \ref{12mono2}.

By slightly modifying the proof of Theorem \ref{12mono2}, we are able to get the monotonicity formula for the  quadharmonic map, i.e.,
\be\nonumber
\Delta^4 u=0.
\ee
Indeed, let $p\rightarrow+\infty$ in \eqref{12Exu} and denote $E_\infty(\la,x,u)=\lim_{p\rightarrow\infty}E(\la,x,u)$,   where the term $\frac{1}{p+1}\la^{8\frac{p+1}{p-1}-n}\int_{\pa B_\la}|u^\la|^{p+1}$ is understood   vanished, then we have
\bc
Assume that $9\leq n\leq17$, then there exist $c_{ij}$ such that  $E_\infty(\la,x,u)$, is a nondecreasing function of $\la>0$. Furthermore,
\be\nonumber
\frac{d E_\infty(\la,x,u)}{d\la}\geq C(n)\la^{-n}
\int_{\pa B_\la(x_0)}\Big(\la\pa_r u\Big)^2,
\ee
where $C(n)>0$ is a constant independent of $\la$.

\ec

\subsection{The calculation of $\frac{d}{d\la}\overline{E}(u,\la)$}

Suppose that $x=0$ in the functional $E(\lambda, x, u)$ and denote by $B_\la$  the balls  centered at zero with radius $\la$. Set
\be\nonumber
\overline{E}(u,\la):=\la^{8\frac{p+1}{p-1}-n}\Big(\int_{B_\la}\frac{1}{2}|\Delta^2 u|^2
-\frac{1}{p+1}\int_{B_\la}|u|^{p+1}\Big).\\
\ee
Set \be
v:=\Delta u, w:=\Delta v, z:=\Delta w.
\ee
Define
\be\label{11b000}\aligned
u^\la(x):&=\la^{\frac{8}{p-1}}(\la x), v^\la(x):=\la^{\frac{8}{p-1}+2}v(\la x),\\
w^\la(x):&=\la^{\frac{8}{p-1}+4} w(\la x), z^\la(x):=\la^{\frac{8}{p-1}+6}z(\la x).
\endaligned\ee
 Therefore,
\be\label{11delta}
\Delta u^\la(x)=v^\la(x), \Delta v^\la(x)=w^\la(x), \Delta w^\la(x)=z^\la.
\ee
Furthermore, differentiating \eqref{11b000} with respect to $\la$ we have
\be\nonumber
\Delta \frac{du^\la}{d\la}=\frac{dv^\la}{d\la}, \Delta \frac{dv^\la}{d\la}=\frac{d w^\la}{d\la}, \Delta \frac{d w^\la}{d\la}=\frac{dz^\la}{d\la}.
\ee
Note that
\be\nonumber
\overline{E}(u,\la)=\overline{E}(u^\la,1)=\int_{B_1}\frac{1}{2}|\Delta^2 u^\la|^2-\frac{1}{p+1}\int_{B_1}|u^\la|^{p+1}.
\ee

Taking derivative of the energy with respect to $\la$ and integrating by part we have
\be\nonumber\aligned
\frac{d\overline{E}(u^\la,1)}{d\la}&=\int_{B_1} \Delta^2 u^\la \Delta^2 \frac{du^\la}{d\la}-\int_{ B_1}|u^\la|^{p-1} u^\la \frac{d u^\la}{d\la}\\
&=\int_{\pa B_1}\Delta^2 u^\la\frac{\pa\Delta\frac{du^\la}{d\la}}{\pa n}-\int_{B_1}\nabla\Delta^2 u^\la\nabla\Delta\frac{du^\la}{d\la}-\int_{ B_1}|u^\la|^{p-1} u^\la \frac{d u^\la}{d\la}\\
&=\int_{\pa B_1}\Delta^2 u^\la\frac{\pa\Delta\frac{du^\la}{d\la}}{\pa n}-\int_{\pa B_1}\frac{\pa \Delta^2 u^\la}{\pa n}\Delta\frac{du^\la}{d\la}+\int_{B_1}\Delta^3 u^\la\Delta\frac{du^\la}{d\la}\\
&\quad-\int_{ B_1}|u^\la|^{p-1} u^\la \frac{d u^\la}{d\la}\\
&=\int_{\pa B_1}\Delta^2 u^\la\frac{\pa\Delta\frac{du^\la}{d\la}}{\pa n}-\int_{\pa B_1}\frac{\pa \Delta^2 u^\la}{\pa n}\Delta\frac{du^\la}{d\la}
+\int_{\pa B_1}\Delta^3 u^\la\frac{\pa\frac{du^\la}{d\la}}{\pa n}\\
&\quad-\int_{B_1}\nabla\Delta u^\la\nabla\frac{du^\la}{d\la}-\int_{ B_1}|u^\la|^{p-1} u^\la \frac{d u^\la}{d\la}\\
&=\int_{\pa B_1}\Delta^2 u^\la\frac{\pa\Delta\frac{du^\la}{d\la}}{\pa n}-\int_{\pa B_1}\frac{\pa \Delta^2 u^\la}{\pa n}\Delta\frac{du^\la}{d\la}
+\int_{\pa B_1}\Delta^3 u^\la\frac{\pa\frac{du^\la}{d\la}}{\pa n}\\
&\quad-\int_{\pa B_1}\frac{\pa \Delta^3u^\la}{\pa n}\frac{du^\la}{d\la}+\int_{B_1}\Delta^4 u^\la\frac{du^\la}{d\la}-\int_{ B_1}|u^\la|^{p-1} u^\la \frac{d u^\la}{d\la}.
\endaligned\ee
In view of \eqref{12LE} and \eqref{11delta}, we have the following

\be\nonumber\aligned
\frac{d\overline{E}(u^\la,1)}{d\la}&=\int_{\pa B_1}
\Delta^3u^\la\frac{\pa}{\pa r}\frac{du^\la}{d\la}+\Delta^2 u^\la \frac{\pa}{\pa r}\frac{d}{d\la}\Delta u^\la\\
&-\frac{\pa}{\pa r}\Delta^3 u^\la \frac{du^\la}{d\la}-\frac{\pa}{\pa r}\Delta^2 u^\la \frac{d}{d\la}\Delta u^\la\\
&=\int_{\pa B_1}  z^\la\frac{\pa}{\pa r}\frac{du^\la}{d\la}+ w^\la\frac{\pa}{\pa r}\frac{dv^\la}{d\la}\\
&\quad-\frac{\pa}{\pa r}z^\la\frac{du^\la}{d\la}-\frac{\pa}{\pa r} w^\la \frac{dv^\la}{d\la}.\\\endaligned\ee
Further, from \eqref{11delta}, let $k:=\frac{8}{p-1}$, we have that
\be\nonumber\aligned
\frac{\pa}{\pa r}\frac{du^\la}{d\la}&=\la\frac{d^2u^\la}{d\la^2}-(k-1)\frac{du^\la}{d\la},\\
\frac{\pa}{\pa r}\frac{dv^\la}{d\la}&=\la\frac{d^2v^\la}{d\la^2}-(k+1)\frac{du^\la}{d\la},\\
\endaligned\ee
and
\be\nonumber\aligned
\frac{\pa}{\pa r}z^\la=\la\frac{dz^\la}{d\la}-(k+6)z^\la,\\
\frac{\pa}{\pa r}w^\la=\la\frac{d w^\la}{d\la}-(k+4)w^\la.
\endaligned\ee
Therefore, we obtain that
\be\label{12ed1ed2}\aligned
\frac{d\overline{E}(u^\la,1)}{d\la}&=\underbrace{\int_{\pa B_1}\la z^\la\frac{d^2 u^\la}{d\la^2}+7z^\la\frac{du^\la}{d\la}-\la \frac{dz^\la}{d\la}\frac{du^\la}{d\la}}\\
&+\underbrace{\int_{\pa B_1}\la  w^\la\frac{d^2v^\la}{d\la^2}+3w^\la\frac{dv^\la}{d\la}-\la \frac{d w^\la}{d\la}\frac{dv^\la}{d\la}}\\
&:=\overline{E}_{d_1}(u^\la,1)+\overline{E}_{d_2}(u^\la,1).
\endaligned\ee

\subsection{The computations of $\frac{\pa^j}{\pa r^j}u^\la$ by $\frac{\pa^i}{\pa \la^i}u^\la, j=1,2,3,4,5,6$}

We start our derivation from the following
\be\label{11original}
\la\frac{du^\la}{d\la}=\frac{8}{p-1}u^\la+r\frac{\pa}{\pa r}u^\la.
\ee
Differentiating \eqref{11original} $j(j=1,2,3,4,5)$ times with respect to $\la$ we have

\be\label{11la1}
\la\frac{d^2 u^\la}{d\la^2}+\frac{du^\la}{d\la}=\frac{8}{p-1}\frac{du^\la}{d\la}+r\frac{\pa}{\pa r}\frac{du^\la}{d\la},\quad\;\;
\ee
\be\label{11la2}
\la\frac{d^3 u^\la}{d\la^3}+2\frac{d^2u^\la}{d\la^2}=\frac{8}{p-1}\frac{d^2u^\la}{d\la^2}+r\frac{\pa}{\pa r}\frac{d^2u^\la}{d\la^2},
\ee

\be\label{11la3}
\la\frac{d^4 u^\la}{d\la^4}+3\frac{d^3u^\la}{d\la^3}=\frac{8}{p-1}\frac{d^3u^\la}{d\la^3}+r\frac{\pa}{\pa r}\frac{d^3u^\la}{d\la^3},
\ee

\be\label{11la4}
\la\frac{d^5 u^\la}{d\la^5}+4\frac{d^4u^\la}{d\la^4}=\frac{8}{p-1}\frac{d^4u^\la}{d\la^4}+r\frac{\pa}{\pa r}\frac{d^4u^\la}{d\la^4},
\ee

\be\label{11la5}
\la\frac{d^6 u^\la}{d\la^6}+5\frac{d^5u^\la}{d\la^5}=\frac{8}{p-1}\frac{d^5u^\la}{d\la^5}+r\frac{\pa}{\pa r}\frac{d^5u^\la}{d\la^5}.
\ee
Differentiating \eqref{11original} $j(j=1,2,3,4,5)$ times with respect to $r$ we have
\be\label{11r1}
\la\frac{\pa}{\pa r}\frac{du^\la}{d\la}=(\frac{8}{p-1}+1)\frac{\pa}{\pa r}u^\la+r\frac{\pa^2}{\pa r^2}u^\la,
\ee

\be\label{11r2}
\la\frac{\pa^2}{\pa r^2}\frac{du^\la}{d\la}=(\frac{8}{p-1}+2)\frac{\pa^2}{\pa r^2}u^\la+r\frac{\pa^3}{\pa r^3}u^\la,
\ee

\be\label{11r3}
\la\frac{\pa^3}{\pa r^3}\frac{du^\la}{d\la}=(\frac{8}{p-1}+3)\frac{\pa^3}{\pa r^3}u^\la+r\frac{\pa^4}{\pa r^4}u^\la,
\ee

\be\label{11r4}
\la\frac{\pa^4}{\pa r^4}\frac{du^\la}{d\la}=(\frac{8}{p-1}+4)\frac{\pa^4}{\pa r^4}u^\la+r\frac{\pa^5}{\pa r^5}u^\la,
\ee

\be\label{11r5}
\la\frac{\pa^5}{\pa r^5}\frac{du^\la}{d\la}=(\frac{8}{p-1}+5)\frac{\pa^5}{\pa r^5}u^\la+r\frac{\pa^6}{\pa r^6}u^\la.
\ee
From \eqref{11original}, on $\pa B_1$, we have
\be\nonumber
\frac{\pa u^\la}{\pa r}=\la \frac{du^\la}{d\la}-\frac{8}{p-1}u^\la.
\ee
Next from \eqref{11la1}, on $\pa B_1$,  we derive that
\be\nonumber
\frac{\pa}{\pa r}\frac{d u^\la}{d\la}=\la\frac{d^2u^\la}{d\la^2}+(1-\frac{8}{p-1})\frac{du^\la}{d\la}.
\ee
From \eqref{11r1}, combining with the two equations above, on $\pa B_1$, we get
\be\label{11r31}\aligned
\frac{\pa^2}{\pa r^2} u^\la&=\la\frac{\pa}{\pa r}\frac{du^\la}{d\la}-(1+\frac{8}{p-1})\frac{\pa}{\pa r} u^\la\\
&=\la^2\frac{d^2 u^\la}{d\la^2}-\la\frac{16}{p-1}\frac{du^\la}{d\la}+(1+\frac{8}{p-1})\frac{8}{p-1}u^\la.
\endaligned\ee
Differentiating \eqref{11la1} with respect to $r$, and combining  with \eqref{11la1} and \eqref{11la2}, we get that
\be\label{11r32}\aligned
\frac{\pa^2}{\pa r^2}\frac{du^\la}{d\la}&=\la\frac{\pa}{\pa r}\frac{d^2 u^\la}{d\la^2}-\frac{8}{p-1}\frac{\pa}{\pa r}\frac{du^\la}{d\la}\\
&=\la^2\frac{d^3 u^\la}{d\la^3}+(2-\frac{16}{p-1})\la\frac{d^2 u^\la}{d\la^2}-(1-\frac{8}{p-1})\frac{8}{p-1}\frac{du^\la}{d\la}.
\endaligned\ee
From \eqref{11r2}, on $\pa B_1$, combining with \eqref{11r31} and \eqref{11r32}, we have
\be\label{10ue7}\aligned
\frac{\pa^3}{\pa r^3}u^\la=&\la\frac{\pa^2}{\pa r^2}\frac{du^\la}{d\la}-(2+\frac{8}{p-1})\frac{\pa^2}{\pa r^2}u^\la\\
=&\la^3\frac{d^3 u^\la}{d\la^3}-\la^2\frac{24}{p-1}\frac{d^2 u^\la}{d\la^2}+\la(\frac{24}{p-1}+\frac{192}{(p-1)^2})\frac{du^\la}{d\la}\\
&-(2+\frac{8}{p-1})(1+\frac{8}{p-1})\frac{8}{p-1}u^\la.
\endaligned\ee
Now differentiating \eqref{11la1} once with respect to $r$, we get
\be\nonumber
\la\frac{\pa^2}{\pa r^2}\frac{d^2 u^\la}{d\la^2}=(\frac{8}{p-1}+1)\frac{\pa^2}{\pa r^2}\frac{du^\la}{d\la}+r\frac{\pa^3}{\pa r^3}\frac{du^\la}{d\la},
\ee
then on $\pa B_1$, we have
\be\label{11ue8}
\frac{\pa^3}{\pa r^3}\frac{du^\la}{d\la}=\la\frac{\pa^2}{\pa r^2}\frac{d^2 u^\la}{d\la^2}-(\frac{8}{p-1}+1)\frac{\pa^2}{\pa r^2}\frac{du^\la}{d\la}.
\ee
Now differentiating \eqref{11la2} twice with respect to $r$, we get
\be\nonumber
\la\frac{\pa}{\pa r}\frac{d^3 u^\la}{d\la^3}=(\frac{8}{p-1}-1)\frac{\pa}{\pa r}\frac{d^2 u^\la}{d\la^2}+r\frac{\pa^2}{\pa r^2}\frac{d^2 u^\la}{d\la^2},
\ee
hence on $\pa B_1$, combining with \eqref{11la2} and \eqref{11la3} there holds
\be\label{11ue9}
\aligned
\frac{\pa^2}{\pa r^2}&\frac{d^2 u^\la}{d\la^2}=\la\frac{\pa}{\pa r}\frac{d^3 u^\la}{d\la^3}+(1-\frac{8}{p-1})\frac{\pa}{\pa r}\frac{d^2 u^\la}{d\la^2}\\
=&\la^2\frac{d^4 u^\la}{d\la^4}+\la(4-\frac{16}{p-1})\frac{d^3 u^\la}{d\la^3}+(1-\frac{8}{p-1})(2-\frac{8}{p-1})\frac{d^2 u^\la}{d\la^2}.
\endaligned\ee
Now differentiating \eqref{11la1} with respect to $r$, we have
\be\nonumber
\la\frac{\pa}{\pa r}\frac{d^2 u^\la}{d\la^2}=\frac{8}{p-1}\frac{\pa}{\pa r}\frac{du^\la}{d\la}+r\frac{\pa^2}{\pa r^2}\frac{du^\la}{d\la}.
\ee
This combining with \eqref{11la1} and \eqref{11la2}, on $\pa B_1$, we have
\be\label{11ue10}\aligned
\frac{\pa^2}{\pa r^2}\frac{du^\la}{d\la}&=\la\frac{\pa}{\pa r}\frac{d^2 u^\la}{d\la^2}-\frac{8}{p-1}\frac{\pa}{\pa r}\frac{d u^\la}{d\la}\\
&=\la^2\frac{d^3 u^\la}{d\la^3}+\la(2-\frac{16}{p-1})\frac{d^2 u^\la}{d\la^2}-\frac{8}{p-1}(1-\frac{8}{p-1})\frac{du^\la}{d\la}.
\endaligned\ee
Now from \eqref{11ue8}, combining with \eqref{11ue9} and \eqref{11ue10}, we get
\be\label{11ue11}\aligned
\frac{\pa^3}{\pa r^3}\frac{du^\la}{d\la}=&\la^3\frac{d^4 u^\la}{d\la^4}+\la^2(3-\frac{24}{p-1})\frac{d^3 u^\la}{d\la^3}-\la(1-\frac{8}{p-1})\frac{24}{p-1}\frac{d^2 u^\la}{d\la^2}\\
&+(1-\frac{8}{p-1})(1+\frac{8}{p-1})\frac{8}{p-1}\frac{du^\la}{d\la}.
\endaligned\ee
From \eqref{11r3}, on $\pa B_1$, combining with \eqref{11ue11}  we obtain that
\be\nonumber\aligned
\frac{\pa^4}{\pa r^4}u^\la=&\la \frac{\pa^3}{\pa r^3}\frac{du^\la}{d\la}-(3+\frac{8}{p-1})\frac{\pa^3}{\pa r^3}u^\la\\
=&\la^4\frac{d^4 u^\la}{d\la^4}-\la^3\frac{32}{p-1}\frac{d^3 u^\la}{d\la^3}+\la^2(2+\frac{16}{p-1})\frac{24}{p-1}\frac{d^2 u^\la}{d\la^2}\\
&-\la (1+\frac{8}{p-1})(1+\frac{4}{p-1})\frac{64}{p-1}\frac{du^\la}{d\la}\\
&+(3+\frac{8}{p-1})(2+\frac{8}{p-1})(1+\frac{8}{p-1})\frac{8}{p-1}u^\la.\\
\endaligned\ee
From \eqref{11r4}, on $\pa B_1$, we have
\be\nonumber
\frac{\pa^5}{\pa r^5}u^\la=\la\frac{\pa^4}{\pa r^4}\frac{du^\la}{d\la}-(\frac{8}{p-1}+4)\frac{\pa^4}{\pa r^4}u^\la.
\ee
Now differentiating \eqref{11la1} three and four times with respect to $r$, we get
\be\label{11la1r3}
\la\frac{\pa^3}{\pa r^3}\frac{d^2 u^\la}{d\la^2}=(\frac{8}{p-1}+2)\frac{\pa^3}{\pa r^3}\frac{du^\la}{d\la}+r\frac{\pa^4}{\pa r^4}\frac{du^\la}{d\la}
\ee
and
\be\label{11la1r4}
\la\frac{\pa^4}{\pa r^4}\frac{d^2u^\la}{d\la^2}=(\frac{8}{p-1}+3)\frac{\pa^4}{\pa r^4}\frac{du^\la}{d\la}+r\frac{\pa^5}{\pa r^5}\frac{du^\la}{d\la}.
\ee
Next differentiating \eqref{11la2} two, three times with respect to $r$, we have
\be\label{11la2r2}
\la\frac{\pa^2}{\pa r^2}\frac{d^3 u^\la}{d\la^3}=\frac{8}{p-1}\frac{\pa^2}{\pa r^2}\frac{d^2u^\la}{d\la^2}+r\frac{\pa^3}{\pa r^3}\frac{d^2u^\la}{d\la^2}
\ee
and
\be\label{11la2r3}
\la\frac{\pa^3}{\pa r^3}\frac{d^3 u^\la}{d\la^3}=(\frac{8}{p-1}+1)\frac{\pa^3}{\pa r^3}\frac{d^2u^\la}{d\la^2}+r\frac{\pa^4}{\pa r^4}\frac{d^2u^\la}{d\la^2}.
\ee
Differentiating \eqref{11la3} with respect to $r$, we have
\be\label{11la3r1}
r\frac{\pa^2}{\pa r^2}\frac{d^3u^\la}{d\la^3}=\la\frac{\pa}{\pa r}\frac{d^4 u^\la}{d\la^4}+(2-\frac{8}{p-1})\frac{\pa}{\pa r}\frac{d^3u^\la}{d\la^3}.
\ee
From \eqref{11la3} and \eqref{11la4}, on $\pa B_1$, we have
\be\nonumber
\frac{\pa}{\pa r}\frac{d^3u^\la}{d\la^3}=\la\frac{d^4u^\la}{d\la^4}+(3-\frac{8}{p-1})\frac{d^3u^\la}{d\la^3},
\ee
\be\nonumber
\frac{\pa}{\pa r}\frac{d^4u^\la}{d\la^4}=\la\frac{d^5u^\la}{d\la^5}+(4-\frac{8}{p-1})\frac{d^4u^\la}{d\la^4}.
\ee
Hence on the boundary $\pa B_1$, we have
\be\aligned
\frac{\pa^2}{\pa r^2}\frac{d^3u^\la}{d\la^3}=&\la^2\frac{d^5 u^\la}{d\la^5}+(4-\frac{8}{p-1})\la\frac{d^4u^\la}{d\la^4}
+(2-\frac{8}{p-1})\la\frac{d^4u^\la}{d\la^4}\\
&\quad+(3-\frac{8}{p-1})(2-\frac{8}{p-1})\frac{d^3u^\la}{d\la^3}.
\endaligned\ee
Combining with \eqref{11la2r2}, on the boundary $\pa B_1$, we have
\be\aligned
\frac{\pa^3}{\pa r^3}\frac{d^2u^\la}{d\la^2}&=\la\frac{\pa^2}{\pa r^2}\frac{d^3 u^\la}{d\la^3}-\frac{8}{p-1}\frac{\pa^2}{\pa r^2}\frac{d^2 u^\la}{d\la^2}\\
&=\la^3\frac{d^5u^\la}{d\la^5}+(6-\frac{8}{p-1})\la^2\frac{d^4u^\la}{d\la^4}+(3(\frac{8}{p-1})^2-\frac{72}{p-1}+6)\la\frac{d^3u^\la}{d\la^3}\\
&\quad-(1-\frac{8}{p-1})(2-\frac{8}{p-1})\frac{8}{p-1}\frac{d^2u^\la}{d\la^2}.
\endaligned\ee
Hence on the boundary $\pa B_1$, we get that
\be\aligned
\frac{\pa^4}{\pa r^4}\frac{du^\la}{d\la}&=\la\frac{\pa^3}{\pa r^3}\frac{d^2 u^\la}{d\la^2}-(\frac{8}{p-1}+2)\frac{\pa^3}{\pa r^3}\frac{du^\la}{d\la}\\
&=\la^4\frac{d^5u^\la}{d\la^5}+(4-\frac{32}{p-1})\la^3\frac{d^4u^\la}{d\la^4}+(5(\frac{8}{p-1})^2-6\frac{8}{p-1})\la^2\frac{d^3u^\la}{d\la^3}\\
&\quad-4\frac{8}{p-1}(\frac{8}{p-1}-1)(\frac{8}{p-1}+1)\la\frac{d^2u^\la}{d\la^2}\\
&\quad-(\frac{8}{p-1}+2)(1-\frac{8}{p-1})(1+\frac{8}{p-1})\frac{8}{p-1}\frac{du^\la}{d\la}.
\endaligned\ee
Therefore, we obtain that
\be\aligned
\frac{\pa^5}{\pa r^5}u_e^\la&=\la\frac{\pa^4}{\pa r^4}\frac{du^\la}{d\la}-(\frac{8}{p-1}+4)\frac{\pa^4}{\pa r^4}u^\la\\
&=\la^5\frac{d^5u^\la}{d\la^5}-5\frac{8}{p-1}\la^4\frac{d^4u^\la}{d\la^4}+10\frac{8}{p-1}(\frac{8}{p-1}+1)\la^3\frac{du^\la}{d\la^3}\\
&\quad-10\frac{8}{p-1}(\frac{8}{p-1}+1)(\frac{8}{p-1}+2)\la^2\frac{d^2u^\la}{d\la^2}\\
&\quad+5\frac{8}{p-1}(\frac{8}{p-1}+1)(\frac{8}{p-1}+2)(\frac{8}{p-1}+3)\la\frac{du^\la}{d\la}\\
&\quad-\frac{8}{p-1}(\frac{8}{p-1}+1)(\frac{8}{p-1}+2)(\frac{8}{p-1}+3)(\frac{8}{p-1}+4)u^\la.
\endaligned\ee
Finally  we compute $\frac{\pa^6}{\pa r^6}u^\la$.
From \eqref{11r5}, on the boundary $\pa B_1$, we have
\be\nonumber
\frac{\pa^6}{\pa r^6}u^\la=\la\frac{\pa^5}{\pa r^5}\frac{du^\la}{d\la}-(\frac{8}{p-1}+5)\frac{\pa^5}{\pa r^5}u^\la.
\ee
Differentiating \eqref{11la1} four times with respect to $r$, on the boundary $\pa B_1$, we have
\be\label{11la1r4}
\frac{\pa^5}{\pa r^5}\frac{du^\la}{d\la}=\la\frac{\pa^4}{\pa r^4}\frac{d^3u^\la}{d\la^2}-(\frac{8}{p-1}+3)\frac{\pa^4}{\pa r^4}\frac{du^\la}{d\la}.
\ee
Differentiating \eqref{11r4} with respect to $r$, on the boundary $\pa B_1$, we have
\be\label{11la4r1}
\frac{\pa^2}{\pa r^2}\frac{d^4u^\la}{d\la^4}=\la\frac{\pa}{\pa r}\frac{d^5u^\la}{d\la^5}+(3-\frac{8}{p-1})\frac{\pa}{\pa r}\frac{d^4u^\la}{d\la^4}.
\ee
In a similar way, differentiating \eqref{11la2} with respect to $r$ three times and differentiating \eqref{11la3} with respect to $r$ two times , on the boundary $\pa B_1$, we have
\be\label{11la2r3}
\frac{\pa^4}{\pa r^4}\frac{du^\la}{d\la^2}=\la\frac{\pa^3}{\pa r^3}\frac{d^3u^\la}{d\la^3}-(\frac{8}{p-1}+1)\frac{\pa^3}{\pa r^3}\frac{d^2u^\la}{d\la^2}
\ee
and
\be\label{11la3r2}
\frac{\pa^3}{\pa r^3}\frac{d^3u^\la}{d\la^3}=\la\frac{\pa^2}{\pa r^2}\frac{d^4u^\la}{d\la^4}+(1-\frac{8}{p-1})\frac{\pa^2}{\pa r^2}\frac{d^3u^\la}{d\la^3}.
\ee
From \eqref{11la4r1}, combining with \eqref{11la4} and \eqref{11la5}, on the boundary $\pa B_1$, we get that
\be\nonumber
\frac{\pa^2}{\pa r^2}\frac{d^4u^\la}{d\la^4}=\la^2\frac{d^6u^\la}{d\la^6}+(8-2\frac{8}{p-1})\la\frac{d^5u^\la}{d\la^5}
+(\frac{8}{p-1}-4)(\frac{8}{p-1}-3)\frac{d^4u^\la}{d\la^4}.
\ee
Therefore, combining with \eqref{11la3r2}, we get
\be\nonumber\aligned
\frac{\pa^3}{\pa r^3}\frac{d^3u^\la}{d\la^3}&=\la^3\frac{d^6u^\la}{d\la^6}+(9-3\frac{8}{p-1})\la^2\frac{d^5u^\la}{d\la^5}
+3(\frac{8}{p-1}-2)(\frac{8}{p-1}-3)\la\frac{d^4u^\la}{d\la^4}\\
&\quad+(3-\frac{8}{p-1})(2-\frac{8}{p-1})(1-\frac{8}{p-1})\frac{d^3u^\la}{d\la^3}.
\endaligned\ee
This above combining with \eqref{11la2r3} yields
\be\nonumber\aligned
\frac{\pa^4}{\pa r^4}\frac{d^2u^\la}{d\la^2}&=\la\frac{\pa^3}{\pa r^3}\frac{d^3u^\la}{d\la^3}-(\frac{8}{p-1}+1)\frac{\pa^3}{\pa r^3}\frac{d^2u^\la}{d\la^2}\\
&=\la^4\frac{d^6u^\la}{d\la^6}+(8-4\frac{8}{p-1})\la^3\frac{d^5u^\la}{d\la^5}+6(\frac{8}{p-1}-2)(\frac{8}{p-1}-1)\la^2\frac{d^4u^\la}{d\la^4}\\
&\quad-4\frac{8}{p-1}(\frac{8}{p-1}-1)(\frac{8}{p-1}-2)\la\frac{d^3u^\la}{d\la^3}\\
&\quad+\frac{8}{p-1}(\frac{8}{p-1}-1)(\frac{8}{p-1}-2)(\frac{8}{p-1}+1)\frac{d^2u^\la}{d\la^2}.
\endaligned\ee
 Combining with \eqref{11la1r4}, we get that
\be\nonumber\aligned
\frac{\pa^5}{\pa r^5}\frac{du^\la}{d\la}&=\la\frac{\pa^4}{\pa r^4}\frac{d^2u^\la}{d\la^2}-(\frac{8}{p-1}+3)\frac{\pa^4}{\pa r^4}\frac{du^\la}{d\la}\\
&=\la^5\frac{d^6u^\la}{d\la^6}+(5-5\frac{8}{p-1})\la^4\frac{d^5u^\la}{d\la^5}+10\frac{8}{p-1}(\frac{8}{p-1}-1)\la^3\frac{d^4u^\la}{d\la^4}\\
&\quad-10\frac{8}{p-1}(\frac{8}{p-1}-1)(\frac{8}{p-1}+1)\la^2\frac{d^3u^\la}{d\la^3}\\
&\quad+5\frac{8}{p-1}(\frac{8}{p-1}-1)(\frac{8}{p-1}+1)(\frac{8}{p-1}+2)\la\frac{d^2u^\la}{d\la^2}\\
&\quad+(\frac{8}{p-1}+3)(\frac{8}{p-1}+2)(1-\frac{8}{p-1})(\frac{8}{p-1}+1)\frac{8}{p-1}\frac{du^\la}{d\la}.
\endaligned\ee
Therefore,  from \eqref{11r5}, we obtain that
\be\nonumber\aligned
\frac{\pa^6}{\pa r^6}u^\la&=\la\frac{\pa^5}{\pa r^5}\frac{du^\la}{d\la}-(\frac{8}{p-1}+5)\frac{\pa^5}{\pa r^5}u^\la\\
&=\la^6\frac{d^6u^\la}{d\la^6}-6\frac{8}{p-1}\la^5\frac{d^5u^\la}{d\la^5}+15\frac{8}{p-1}(\frac{8}{p-1}+1)\la^4\frac{d^4u^\la}{d\la^4}\\
&\quad-20\frac{8}{p-1}(\frac{8}{p-1}+1)(\frac{8}{p-1}+2)\la^3\frac{d^3u^\la}{d\la^3}\\
&\quad+15\frac{8}{p-1}(\frac{8}{p-1}+1)(\frac{8}{p-1}+2)(\frac{8}{p-1}+3)\la^2\frac{d^2u^\la}{d\la^2}\\
&\quad-6\frac{8}{p-1}(\frac{8}{p-1}+1)(\frac{8}{p-1}+2)(\frac{8}{p-1}+3)(\frac{8}{p-1}+4)\la\frac{du^\la}{d\la}\\
&\quad+\frac{8}{p-1}(\frac{8}{p-1}+1)(\frac{8}{p-1}+2)(\frac{8}{p-1}+3)(\frac{8}{p-1}+4)(\frac{8}{p-1}+5)u^\la.
\endaligned\ee
In a summary, let $k:=\frac{8}{p-1}$, we have the following
\be\label{12par12}\aligned
\frac{\pa u^\la}{\pa r}&=\la\frac{du^\la}{d\la}-ku^\la,\\
\frac{\pa^2 u^\la}{\pa r^2}&=\la^2\frac{d^2u^\la}{d\la^2}-2k\la\frac{du^\la}{d\la}+k(k+1)u^\la,\\
\endaligned\ee

\be\label{12par34}\aligned
\frac{\pa^3 u^\la}{\pa r^3}&=\la^3\frac{d^3u^\la}{d\la^3}-3k\la^2\frac{d^3u^\la}{d\la^2}+3k(k+1)\la\frac{du^\la}{d\la}-k(k+1)(k+2)u^\la,\\
\frac{\pa^4 u^\la}{\pa r^4}&=\la^4\frac{d^4u^\la}{d\la^4}-4k\la^3\frac{d^3u^\la}{d\la^3}+6k(k+1)\la^2\frac{d^2 u^\la}{d\la^2}\\
&\quad-4k(k+1)(k+2)\la\frac{d u^\la}{d\la}+k(k+1)(k+2)(k+3)u^\la,\\
\endaligned\ee

\be\label{12par5}\aligned
\frac{\pa^5 u^\la}{\pa r^5}&=\la^5\frac{d^5u^\la}{d\la^5}-5k\la^4\frac{d^4u^\la}{d\la^4}+10k(k+1)\la^3\frac{d^3u^\la}{d\la^3}
-10k(k+1)(k+2)\la^2\frac{d^2u^\la}{d\la^2} \\
&\quad+5k(k+1)(k+2)(k+3)\la\frac{du^\la}{d\la}\\
&\quad-k(k+1)(k+2)(k+3)(k+4)u^\la,\\
\endaligned\ee
\be\label{12par6}\aligned
\frac{\pa^6u^\la}{\pa r^6}&=\la^6\frac{d^6u^\la}{d\la^6}-6k\la^5\frac{d^5u^\la}{d\la^5}+15k(k+1)\la^4\frac{d^4u^\la}{d\la^4}
-20k(k+1)(k+2)\la^3\frac{d^3u^\la}{d\la^3}\\
&\quad+15k(k+1)(k+2)(k+3)\la^2\frac{d^2u^\la}{d\la^2}\\
&\quad-6k(k+1)(k+2)(k+3)(k+4)\la\frac{du^\la}{d\la}\\
&\quad+k(k+1)(k+2)(k+3)(k+4)(k+5)u^\la.
\endaligned\ee

\section{On the operators $\Delta^j,j=1,2,3$}

 Let us recall that
 \be\nonumber\aligned
 \Delta u=(\pa_{rr}+\frac{n-1}{r}\pa_r)u+\Delta_\theta (r^{-2}u),
 \endaligned\ee

 \be\nonumber\aligned
\Delta^2 u=&(\pa_{rr}+\frac{n-1}{r}\pa_r)^2+\Delta_\theta\Big(
(\pa_{rr}+\frac{n-1}{r}\pa_r)(r^{-2}u)+r^{-2}(\pa_{rr}+\frac{n-1}{r}\pa_r)u\Big)\\
&\;+\Delta^2_\theta(r^{-4}u)
\endaligned\ee
and
  \be\label{12ddd3}\aligned
\Delta^3u&=(\pa_{rr}+\frac{n-1}{r}\pa_r)^3+\Delta_\theta\Big(
(\pa_{rr}+\frac{n-1}{r}\pa_r)\big(r^{-2}(\pa_{rr}+\frac{n-1}{r}\pa_r)u\big)\\
&\quad+(\pa_{rr}+\frac{n-1}{r}\pa_r)^2(r^{-2}u)+
r^{-2}(\pa_{rr}+\frac{n-1}{r}\pa_r)^2
\Big)
+\Delta^2_\theta\Big(
(\pa_{rr}\\
&\quad +\frac{n-1}{r}\pa_r)(r^{-4}u)+r^{-4}(\pa_{rr}+\frac{n-1}{r}\pa_r)u+r^{-2}(\pa_{rr}+\frac{n-1}{r}\pa_r)(r^{-2}u)
\Big)\\
&\quad+\Delta^3_\theta(r^{-6}u)\\
&:=F_0(u)+\Delta_\theta F_1(u)+\Delta_\theta^2 F_2(u)+\Delta_\theta^3F_3(u).
 \endaligned\ee
 By a direct calculation, if we denote that $a=n-1$, on the boundary $\pa B_1$, then we have
 \be\nonumber\aligned
F_0(u):=&(\pa_{rr}+\frac{n-1}{r}\pa_r)^3u=\Big(\pa_{r^6}+3a\pa_{r^5}+3a(a-2)\pa_{r^4}+a(a-2)(a-7)\pa_{r^3}\\
&-3a(a-2)(a-4)\pa_{r^2}+3a(a-2)(a-4)\pa_r\Big)u\\
:=&\sum_{j=1}^6a_j\pa_{r^j}u,
\endaligned\ee
 \be\nonumber\aligned
F_1(u):=&\Big(3\pa_{r^4}+(6a-12)\pa_{r^3}+(3a^2-24a+42)\pa_{r^2}+(60a-9a^2-96)\pa_r+\\
&(8a^2-64a+120)\Big)u\\
:=&\sum_{j=0}^4b_j\pa_{r^j}u,
 \endaligned\ee
 and
  \be\nonumber\aligned
F_3(u):=&\Big(3\pa_{r^2}+(3a-12)\pa_r+26-6a
\Big)u=\sum_{j=0}^2v_j\pa_{r^j},
 \endaligned\ee
here $\pa_{r^3}:=\pa_{rrr}$ and so on.

\vskip0.13in
 Let us recall that \eqref{12par12}, \eqref{12par34}, \eqref{12par5} and \eqref{12par6} in the end of the previous section, we can turn the differential
with  respect to $r$ into with  respect to $\la$. Thus we have the following
 \be\label{12f012}\aligned
 F_0(u)=\sum_{j=0}^6k_j \la^j\frac{d^ju^\la}{d\la^j},F_1(u)=\sum_{j=0}^4(-t_j)\la^j\frac{d^ju^\la}{d\la^j},F_2(u)=\sum_{j=0}^2
 e_j\la^j\frac{d^ju^\la}{d\la^j}.
  \endaligned\ee
For the simplicity, we let
 \be\nonumber\aligned
a_6&=1,a_5=3a,a_4=3a(a-2),a_3=a(a-2)(a-7),a_2=-3a(a-2)(a-4),\\
a_1&=3a(a-2)(a-4).
  \endaligned\ee
Then we have $k_j$ determined by
 \be\label{12kj}\aligned
k_6&=1,k_5=-6ka_6+a_5,k_4=15k(k+1)a_6-5ka_5+a_4,\\
k_3&=-20k(k+1)(k+2)a_6+10k(k+1)a_5-4ka_4+a_3,\\
k_2&=15k(k+1)(k+2)(k+3)a_6-10k(k+1)(k+2)a_5+6k(k+1)a_4-3ka_3+a_2,\\
k_1&=-6k(k+1)(k+2)(k+3)(k+4)a_6+5k(k+1)(k+2)(k+3)a_5\\
&\quad-4k(k+1)(k+2)a_4+3k(k+1)a_3-2ka_2+a_1,\\
k_0&=k(k+1)(k+2)(k+3)(k+4)(k+5)a_6-k(k+1)(k+2)(k+3)(k+4)a_5\\
&\quad+k(k+1)(k+2)(k+3)a_4-k(k+1)(k+2)a_3+k(k+1)a_2-k a_1;
  \endaligned\ee
and $t_j$ are  determined by
  \be\label{12tj}\aligned
t_4&=-b_4, t_3=4b_4k-b_3,t_2=-6b_4k(k+1)+3b_3k-b_2,\\
t_1&=4b_4k(k+1)(k+2)-3b_3k(k+1)+2b_2k-b_1,\\
t_0&=-b_4k(k+1)(k+2)(k+3)+b_3k(k+1)(k+2)-b_2k(k+1)+b_1k-b_0
  \endaligned\ee
and $e_j$ are determined by
 \be\label{12ej}\aligned
 e_2=3,e_1=-6k+3a-12,e_0= 3k(k+1)-(3a-12)k+26-6a.
 \endaligned\ee


\section{Differentiating  by part formulas}
In all this sections, we denote that $f^{(j)}=\frac{d^j f}{d\la^j}$ and respectively.

\bl\label{12fd1} We have the following type-1 (i.e., $\la^jf^{(j)}f^{(1)}$) differentiating by part formulas:
\be\nonumber\aligned
f f^{(1)}&=\frac{d}{d\la}(\frac{1}{2}f^2),\\
\la^2f^{(2)}f^{(1)}&=-\la(f^{(1)})^2+\frac{d}{d\la}(\frac{1}{2}\la^2f^{(1)}f^{(1)}),\\
\la^3f^{(3)}f^{(1)}&=3\la(f^{(1)})^2-\la^3(f^{(2)})^2+\frac{d}{d\la}(\la^3f^{(2)}f^{(1)}),\\
\endaligned\ee
 \be\nonumber\aligned
\la^4f^{(4)}f^{(1)}=&-12\la(f^{(1)})^2+6\la^3(f^{(2)})^2+\frac{d}{d\la}\big(
\la^4f^{(3)}f^{(1)}-\frac{1}{2}\la^4f^{(2)}f^{(2)}\\
&\;\;-4\la^3f^{(2)}f^{(1)}+6\la^2f^{(1)}f^{(1)}
\big),\\
\endaligned\ee

  \be\nonumber\aligned
\la^5f^{(5)}f^{(1)}=&60\la(f^{(1)})^2-40\la^3(f^{(2)})^2+\la^5(f^{(3)})^2
+\frac{d}{d\la}\big(
\la^5f^{(4)}f^{(1)}-\la^5f^{(3)}f^{(2)}\\
&\;\;-5\la^4f^{(3)}f^{(1)}+5\la^4f^{(2)}f^{(2)}+20\la^3f^{(2)}f^{(1)}-30\la^2f^{(1)}f^{(1)}
\big),
\endaligned\ee
  \be\nonumber\aligned
\la^6f^{(6)}f^{(1)}=&-360\la(f^{(1)})^2+300\la^3(f^{(2)})^2-14\la^5(f^{(3)})^2
+\frac{d}{d\la}\Big(
\la^6f^{(5)}f^{(1)}\\
&\;\;-6\la^5f^{(4)}f^{(1)}+12\la^5f^{(3)}f^{(2)}+30\la^4f^{(3)}f^{(1)}-45\la^4f^{(2)}f^{(2)}\\
&\;\;-120\la^3f^{(2)}f^{(1)}+180\la^2f^{(1)}f^{(1)}-\la^6f^{(4)}f^{(2)}+\frac{1}{2}\la^6f^{(3)}f^{(3)}
\Big),
\endaligned\ee
  \be\nonumber\aligned
\la^7f^{(7)}f^{(1)}=&2520\la(f^{(1)})^2-2520\la^3(f^{2})^2+189\la^5(f^{(3)})^2-\la^7(f^{(4)})^2\\
&+\frac{d}{d\la}\Big(
\la^7f^{(6)}f^{(1)}-7\la^6f^{(5)}f^{(1)}+42\la^5f^{(4)}f^{(1)}-84\la^5f^{(3)}f^{(2)}\\
&\;\;-210\la^4f^{(3)}f^{(1)}+315\la^4f^{(2)}f^{(2)}+840\la^3f^{(2)}f^{(1)}\\
&\;\;-1260\la^2f^{(1)}f^{(1)}
+7\la^6f^{(4)}f^{(2)}-7\la^6f^{(3)}f^{(3)}\\
&\;\;-\la^7f^{(5)}f^{(2)}+\la^7f^{(4)}f^{(3)}
\Big).
\endaligned\ee
\el
\bp These formulas above can be checked directly. \ep
\br
We see that we can decompose the term $\la^{j}f^{(j)}f^{(1)}$ into two parts, the quadratic form and derivative term, i.e.,

\be\nonumber\aligned
\la^{j}f^{(j)}f^{(1)}=\sum_{s\leq \frac{j+1}{2},s\in N}b_{j,s}\la^{2s-1}(f^{(s)})^2+\frac{d}{d\la}\Big(\sum_{i,l}c_{i,l}\la^{i+l}f^{(i)}f^{(l)}\Big).
\endaligned\ee

\er

\bl\label{12fd2} We have the following type-2 (i.e., $\la^{j+1}f^{(j)}f^{(2)}$) differential by part formulas:
\be\nonumber\aligned
&\la ff^{(2)}=-\la(f^{(1)})^2+\frac{d}{d\la}\big(
\la ff^{(1)}-\frac{1}{2}f^2
\big),\\
&\la^2 f^{(1)}f^{(2)}=-\la(f^{(1)})^2+\frac{d}{d\la}(\frac{1}{2}\la^2f^{(1)}f^{(1)}),\\
&\la^4 f^{(3)}f^{(2)}=-2\la^3(f^{(2)})^2+\frac{d}{d\la}(\frac{1}{2}\la^4f^{(2)}f^{(2)}),\\
&\la^5f^{(4)}f^{(2)}=10\la^3(f^{(2)})^2-\la^5(f^{(3)})^2+\frac{d}{d\la}\big(
\la^5f^{(3)}f^{(2)}-\frac{5}{2}\la^4f^{(2)}f^{(2)}
\big),\\
&\la^6f^{(5)}f^{(2)}=-60\la^3(f^{(2)})^2+9\la^5(f^{(3)})^2
+\frac{d}{d\la}\Big(
\la^6f^{(4)}f^{(2)}-6\la^5f^{(3)}f^{(2)}\\
&\;\;+15\la^4f^{(2)}f^{(2)}-\frac{1}{2}\la^6f^{(3)}f^{(3)}
\Big),\\
&\la^7f^{(6)}f^{(2)}=420\la^3(f^{(2)})^2-84\la^5(f^{(3)})^2+\la^7(f^{(4)})^2
+\frac{d}{d\la}\Big(
\la^7f^{(5)}f^{(2)}\\
&\;\;-\la^7f^{(4)}f^{(3)}+\frac{7}{2}\la^6f^{(3)}f^{(3)}-7\la^6f^{(4)}f^{(2)}
+42\la^5f^{(3)}f^{(2)}\\
&\;\;-105\la^4f^{(2)}f^{(2)}+\frac{7}{2}\la^6f^{(3)}f^{(3)}
\Big).
\endaligned\ee
\el

\bp
These formulas can be verified  directly.
\ep

\br
We note that we can decompose the term $\la^{j+1}f^{(j)}f^{(2)}$ into two parts, the quadratic form and derivative term, i.e.,

\be\nonumber\aligned
\la^{j+1}f^{(j)}f^{(2)}=\sum_{s\leq \frac{j+2}{2},s\in N}a_{j,s}\la^{2s-1}(f^{(s)})^2+\frac{d}{d\la}\Big(\sum_{i,l}c_{i,l}\la^{i+l}f^{(i)}f^{(l)}\Big).
\endaligned\ee

\er

\section{The term $\overline{E}_{d_2}(u^\la,1)$}

Recall the definition in \eqref{12ed1ed2}
\be\label{12ed2}
\overline{E}_{d_2}(u^\la,1)=\int_{\pa B_1}(\la w^\la\frac{d^2v^\la}{d\la^2}+3w^\la\frac{dv^\la}{d\la}-\la\frac{dw^\la}{d\la}\frac{dv^\la}{d\la}).
\ee
Since $w^\la=\Delta v^\la$
$=\pa_{rr}v^\la+\frac{n-1}{r}\pa_r v^\la$
$ +r^{-2}\md(\ds v^\la)$, in view of \eqref{12par12}, on the boundary $\pa B_1$, we have
\be\nonumber\aligned
w^\la&=\la^2\frac{d^2v^\la}{d\la^2}+\la\frac{dv^\la}{d\la}(n-1-\frac{16}{p-1})\\
&\quad+u^\la\frac{8}{p-1}(2+\frac{8}{p-1}-n)+\md(\ds v^\la)\\
&:=\la^2\frac{d^2 v^\la}{d\la^2}+\alpha\la\frac{dv^\la}{d\la}+\beta v^\la+\md(\ds v^\la).\\
\endaligned\ee
Integrate by part with suitable times we have the following
\be\label{12ed2}\aligned
\overline{E}_{d_2}(u^\la,1)&=\int_{\pa B_1}\big[ 2\la^3(\frac{d^2v^\la}{d\la^2})^2+(2\alpha-2\beta-4)\la(\frac{dv^\la}{d\la})^2\big]\\
&\quad\frac{d}{d\la}\int_{\pa B_1}
\Big[\frac{\beta}{2}\frac{d}{d\la}(\la(v^\la)^2)^2-\frac{1}{2}\la^3\frac{d}{d\la}(\frac{dv^\la}{d\la})^2+(\frac{\beta}{2}+2)(v^\la)^2\\
&\quad-\frac{1}{2}\frac{d}{d\la}(\la|\nabla_{\q} v^\la|^2)-\frac{1}{2}|\nabla_{\q} v^\la|^2\Big]\\
&\quad+2\int_{\pa B_1}\la|\nabla_{\q}
\frac{dv^\la}{d\la}|^2.
\endaligned\ee
 Let us further investigate the inner structure of $\overline{E}_{d_2}(u^\la,1)$, we can obtain more and  crucial information for our construction
of the  monotonicity formula under the desired condition.
 Since $v^\la=\Delta u^\la$, on the boundary $\pa B_1$, by a direct calculation we have the following
 \be\nonumber\aligned
\frac{dv^\la}{d\la}=\la^2\frac{d^3u^\la}{d\la^3}+(\alpha+2)\la\frac{d^2u^\la}{d\la^2}+(\alpha+\beta)\frac{du^\la}{d\la}+\Delta_\theta\frac{du^\la}{d\la},
\endaligned\ee
 \be\nonumber\aligned
\frac{d^2v^\la}{d\la^2}=\la^2\frac{d^4u^\la}{d\la^4}+(\alpha+4)\la\frac{d^3u^\la}{d\la^3}+(2\alpha+\beta+2)\frac{d^2u^\la}{d\la^2}+\Delta_\theta\frac{d^2u^\la}{d\la^2},
\endaligned\ee
therefore differential by part, we can get that
 \be\label{12abc}\aligned
&\int_{\pa B_1}(2\alpha-2\beta-4)\la(\frac{dv^\la}{d\la})^2+2\la^3(\frac{d^2v^\la}{d\la^2})^2
=\int_{\pa B_1}\big(2\alpha^2-2\alpha-6\beta-28\big)\la^5(\frac{d^3u^\la}{d\la^3})^2\\
&\quad+\big(2\alpha^3+(-16-2\beta)\alpha^2+16\alpha+6\beta^2+32\beta+40\big)\la^3(\frac{d^2u^\la}{d\la^2})^2+2\la^7(\frac{d^4u^\la}{d\la^4})^2\\
&\quad+\big(-2\alpha^3+(2\beta+8)\alpha^2+(2\beta^2-8)\alpha-2\beta^3-8\beta^2-8\beta\big)\la(\frac{du^\la}{d\la})^2\\
&\quad+\frac{d}{d\la}\Big(\sum_{i,j}c_{i,j}\la^{i+j}\frac{d^iu^\la}{d\la^i}\frac{d^ju^\la}{d\la^j}\Big)+(2\alpha-2\beta-4)\la(\Delta_\theta\frac{du^\la}{d\la})^2\\
&\quad+2\la^5(\nabla_\theta\frac{d^3u^\la}{d\la^3})^2
+(14-2\beta)\la^3(\nabla_\theta\frac{d^2u^\la}{d\la^2})^2
+(-2-2\beta)\la (\nabla_\theta\frac{du^\la}{d\la})^2\\
&\quad+\frac{d}{d\la}\Big(\sum_{i,j}e_{i,j}\la^{i+j}\nabla_\theta\frac{d^iu^\la}{d\la^i}\nabla_\theta\frac{d^ju^\la}{d\la^j}\Big)
+2\la^3(\Delta_\theta\frac{d^2u^\la}{d\la^2})^2\\
&:=\int_{\pa B_1}\sum_{j=1}^4a_j\la^{2j-1}(\frac{d^ju^\la}{d\la^j})^2+\Big(\sum_{s=1}^3b_s\la^{2s-1}(\nabla_\theta\frac{d^su^\la}{d\la^s})^2\Big)
+\sum_{l=1}^2c_l\la^{2l-1}(\Delta_\theta\frac{d^lu^\la}{d\la^l})^2
\\&\quad+\frac{d}{d\la}\Big(\sum_{i,j}c_{i,j}\la^{i+j}\frac{d^iu^\la}{d\la^i}\frac{d^ju^\la}{d\la^j}\Big)\\
&\quad+\frac{d}{d\la}\Big(\sum_{i,j}e_{i,j}\la^{i+j}\nabla_\theta\frac{d^iu^\la}{d\la^i}\nabla_\theta\frac{d^ju^\la}{d\la^j}\Big),\\
\endaligned\ee
where $c_{i,j}$ and $e_{i,j}$ may  have   exact expressions of $\alpha,\beta$, but we do not intend to  give them here since the key term is the quadratic form.

\br From \eqref{12ed2}, the first term of the above integral is positive. Recall that $v^\la=\Delta u^\la$
now we have
\be\aligned
\overline{E}_{d_2}(u^\la,1)&\geq\frac{d}{d\la}\int_{\pa B_1}
\Big[\frac{\beta}{2}\frac{d}{d\la}(\la(\Delta u^\la)^2)^2-\frac{1}{2}\la^3\frac{d}{d\la}(\frac{d\Delta u^\la}{d\la})^2+\frac{\beta}{2}(\Delta_b u^\la)^2\\
&\quad-\frac{1}{2}\frac{d}{d\la}(\la|\nabla_{\q} \Delta u^\la|^2)-\frac{1}{2}|\nabla_{\q} \Delta u^\la|^2\Big].\\
\endaligned\ee
\er
If we use this estimate alone,  we can not construct the desired monotonicity formula for all  $n$ with  $\frac{n+8}{n-8}<p<p_c(n)$. More precisely, when $n\in[15,27]$, it seems that, under the condition $\frac{n+8}{n-8}<p<p_c(n)$,  the desired monotonicity formula can not hold.


\section{ The term $\overline{E}_{d_1}(u^\la,1)$ }

Recall that \eqref{12ed1ed2}, we have
\be\nonumber
\overline{E}_{d_1}(u^\la,1)=\int_{\pa B_1}\big(\la z^\la\frac{d^2u^\la}{d\la^2}+7z^\la\frac{du^\la}{d\la}-\la\frac{dz^\la}{d\la}\frac{du^\la}{d\la}\big).
\ee
\subsection{The integral corresponding to  the operator $F_0$}

First, we consider the operator $F_0$ as defined  in \eqref{12f012}.
We split the integral into two parts, we denote that
\be\nonumber
F_{01}:=\int_{\pa B_1}\big(\la F_0(u^\la)\frac{d^2u^\la}{d\la^2}\big),
\ee
\be\nonumber
F_{02}:=\int_{\pa B_1}\big(7F_0(u^\la)-\la\frac{d F_0(u^\la)}{d\la}\big)\frac{du^\la}{d\la}.
\ee
Recall that \eqref{12f012}, if we denote that $f=u^\la,f'=\frac{du^\la}{d\la}$, we have
\be\nonumber
F_0(u^\la)=k_6\la^6\frac{d^6u^\la}{d\la^6}+k_5\frac{d^5u^\la}{d\la^5}+k_4\la^4\frac{d^4u^\la}{d\la^4}
+k_3\frac{d^3u^\la}{d\la^3}+k_2\frac{d^2u^\la}{d\la^2}+k_1\la\frac{du^\la}{d\la}+k_0u\la.
\ee
Hence,
\be\nonumber\aligned
7F_0(u^\la)-\la\frac{d F_0(u^\la)}{d\la}&=-k_6\la^7\frac{d^7u^\la}{d\la^7}+(k_6-k_5)\la^6\frac{d^6u^\la}{d\la^6}+(2k_5-k_4)\la^5\frac{d^5u^\la}{d\la^5}\\
&\quad+(3k_4-k_3)\la^4\frac{d^4u^\la}{d\la^4}+(4k_3-k_2)\la^3\frac{d^3u^\la}{d\la^3}\\
&\quad+(5k_2-k_1)\la^2\frac{d^2u^\la}{d\la^2}
+(6k_1-k_0)\la\frac{du^\la}{d\la}+7k_0u^\la.
\endaligned\ee
By the following differential identity, combining with  the  differential by part formulas of  Section 4, we have
\be\nonumber\aligned
\la^7& f''''''f''+k_5\la^6f'''''f''+k_4\la^5f''''f''+k_3\la^4f'''f''\\
&\quad +k_2\la^3f''f''+k_1\la^2f'f''+k_0\la ff''\\
&=\Big[\la^7f''''''f''-\la^7f''''f'''+(k_5-7)\la^6f''''f''+(7-\frac{k_5}{2})\la^6(f''')^2\\
&\quad+(-6k_5+k_4+42)\la^5f'''f''+(15k_5-\frac{5}{2}k_4+\frac{1}{2}k_3-105)\la^4(f'')^2\\
&\quad+k_0\la f f'-\frac{1}{2}k_0f^2\Big]'+\la^7(f'''')^2+(9k_5-k_4-84)\la^5(f''')^2\\
&\quad+(-60k_5+10k_4-2k_3+k_2+420)\la^3(f'')^2-k_0\la(f')^2+k_1\la^2f'f''
\endaligned\ee
and
\be\nonumber\aligned
&-\la^7f'''''''f'+(1-k_5)\la^6f''''''f'+(2k_5-k_4)\la^5f'''''f'+(3k_4-k_3)\la^4f''''f'\\
&+(4k_3-k_2)\la^3f'''f'+(5k_2-k_1)\la^2f''f'+(6k_1-k_0)\la f'f'+7k_0ff'\\
&=\Big[-\la^7f''''''f'+\la^7f'''''f''-\la^7f''''f'''+(8-k_5)\la^6f'''''f'+(k_5-15)\la^6 f''''f''\\
&\quad-(14k_5-k_4-138)\la^5f'''f''+(55k_5-\frac{13}{2}k_4+\frac{1}{2}k_3-480)\la^4(f'')^2\\
&\quad+(8k_5-k_4-48)\la^5f''''f'+(-40k_5+8k_4-k_3+240)\la^4f'''f'\\
&\quad+(160k_5-32k_4+8k_3-k_2-960)\la^3f''f'+\frac{7}{2}k_0f^2\Big]'\\
&\quad+\la^7(f'''')^2+(14k_5-k_4-138)\la^5(f''')^2+(22-k_5)\la^6f''''f'''\\
&\quad+(-380k_5+58k_4-10k_3+k_2+2820)\la^3(f'')^2+(6k_1-k_0)\la(f')^2\\
&\quad+(-480k_5+96k_4-24k_3+8k_2-k_1+2880)\la^2f''f'.
\endaligned\ee
Then by the above two identities ($f=u^\la, f'=\frac{d}{d\la}u^\la$), we get   the following integral corresponding to the operator $A$:
\be\nonumber\aligned
\mathcal{F}_0&=\int_{\pa B_1}\Big(\la F_0(u^\la)\frac{d^2u^\la}{d\la^2}
+7F_0(u^\la)\frac{du^\la}{d\la}-\la\frac{d F_0(u^\la)}{d\la}\frac{du^\la}{d\la}\Big)\\
&=\mathcal{F}_{01}+\int_{\pa B_1}\Big[ A_4\la^7(\frac{d^4u^\la}{d\la^4})^2+A_3\la^5(\frac{d^3u^\la}{d\la^3})^2+A_2\la^3(\frac{d^2u^\la}{d\la^2})^2+A_1\la(\frac{du^\la}{d\la})^2\Big ],
\endaligned\ee
where
\be\label{12d1d2d3}\aligned
A_4&=2,\\
A_3&=26k_5-2k_4-288k_6,\\
A_2&=-440k_5+68k_4-12k_3+2k_2+3240k_6,\\
A_1&=480k_5-96k_4+24k_3-8k_2+6k_1-2k_0-2280k_6,\\
\endaligned\ee
and the part $\mathcal{F}_{01}$ denotes the differential term, exactly, it is
\be\nonumber\aligned
\mathcal{F}_{01}:&=\frac{d}{d\la}\int_{\pa B_1}
\Big[-\la^7\frac{d^6u^\la}{d\la^6}\frac{du^\la}{d\la}+2\la^7\frac{d^5u^\la}{d\la^5}\frac{d^2u^\la}{d\la^2}
-2\la^7\frac{d^4u^\la}{d\la^4}\frac{d^3u^\la}{d\la^3}\\
&\quad+(8-k_5)\la^6\frac{d^5u^\la}{d\la^5}\frac{du^\la}{d\la}
+(-20k_5+2k_4+180)\la^5\frac{d^3u^\la}{d\la^3}\frac{d^2u^\la}{d\la^2}\\
&\quad+(70k_5-9k_4+k_3-585)\la^4(\frac{d^2u^\la}{d\la^2})^2+(8k_5-k_4-48)\la^5\frac{d^4u^\la}{d\la^4}\frac{du^\la}{d\la}\\
&\quad+(-40k_5+8k_4-k_3+240)\la^4\frac{d^3u^\la}{d\la^3}\frac{du^\la}{d\la}\\
&\quad+(160k_5-32k_4+8k_3-k_2-960)\la^3\frac{d^2u^\la}{d\la^2}\frac{du^\la}{d\la}+3k_0(u^\la)^2\\
&\quad+(2k_5-22)\la^6\frac{d^4u^\la}{d\la^4}\frac{d^2u^\la}{d\la^2}+(18-k_5)\la^6(\frac{d^3u^\la}{d\la^3})^2+k_0\la u_e^\la\frac{du^\la}{d\la}\\
&\quad+(-240k_5+48k_4-12k_3+4k_2+\frac{1}{2}k_1+1440)\la^2(\frac{du^\la}{d\la})^2\Big].
\endaligned\ee

\subsection{The integrals corresponding to the  operator $\Delta_\theta F_1(u^\la)$}

Let us define

\be\nonumber\aligned
I_1(F_1):=&\int_{\pa B_1}\la^{1} \Delta_\theta F_1(u^\la)\frac{d^2u^\la}{d\la^2}=\int_{\pa B_1}\la^{j+1} \sum_{j=0}^4(-t_j)\Delta_\theta\frac{d^j u^\la}{d\la^j}\frac{d^2u^\la}{d\la^2}\\
=&\int_{\pa B_1}\sum_{j=0}^4 t_j \la^{j+1}\nabla_\theta \frac{d^j u^\la}{d\la^j}\nabla_\theta\frac{d^2u^\la}{d\la^2},
\endaligned\ee
and
\be\nonumber\aligned
I_2(F_1):&=\int_{\pa B_1} \Big(7\Delta_\theta F_1(u^\la)-\la\frac{d}{d\la}\Delta_\theta F_1(u)\Big)\frac{du^\la}{d\la}\\
&=\int_{\pa B_1}\Big(7\Delta_\theta \big(\sum_{j=0}^4(-t_j)\la^j\frac{d^j u^\la}{d\la^j}\big)
-\la\frac{d}{d\la}\Delta_\theta \big(\sum_{j=0}^4(-t_j)\la^j\frac{d^j u^\la}{d\la^j}\big)\Big)\frac{du^\la}{d\la}\\
&=\int_{\pa B_1}\sum_{j=1}^5t_{0j}\la^j\nabla_\theta\frac{d^ju^\la}{d\la^j}\nabla_\theta\frac{du^\la}{d\la},
\endaligned\ee
where
\be\nonumber\aligned
t_{05}=-t_4,t_{04}=3t_4-t_3,t_{03}=4t_3-t_2,t_{02}=5t_2-t_1,t_{01}=6t_1-t_0,
\endaligned\ee
and $t_j$ is  defined in \eqref{12tj}.
Recall the    formulas in the Lemma \ref{12fd1}-\ref{12fd2}, we regard that $f=\nabla_\theta u^\la$, then we can obtain that

\be\nonumber\aligned
I_1(F_1)+I_2(F_1)=&\int_{\pa B_1}B_1\la(\nabla_\theta \frac{du^\la}{d\la})^2+B_2\la^3(\nabla_\theta\frac{d^2u^\la}{d\la^2})^2
+B_3\la^5(\nabla_\theta\frac{d^3u^\la}{d\la^3})^2\\
&\;\;+\frac{d}{d\la}\int_{\pa B_1}\Big(\sum_{0\leq i,j\leq 2}b_{i,j}\la^{i+j}\nabla_\theta\frac{d^iu^\la}{d\la^i}\nabla_\theta\frac{d^ju^\la}{d\la^j}\Big),
\endaligned\ee
where
\be\label{12B123}\aligned
B_1=-2t_0+6t_1-8t_2+24t_3-96t_4,
B_2=2t_2-12t_3+68t_4,
B_3=-2t_4=6.
\endaligned\ee
The coefficient $b_{i,j}$ can be determined by $t_j$ but we do not give the precise form since the constant  is not important for our estimate below.

\subsection{The integral corresponding to the operator $\Delta_\theta^2F_2(u^\la)$}

Recall that the sphere representation of triple-harmonic operator, i.e., \eqref{12ddd3} and \eqref{12f012}.
Let us define the following
\be\nonumber\aligned
I_1(F_2):&=\int_{\pa B_1}\la^{j+1} \Delta^2_\theta F_2(u^\la)\frac{d^2u^\la}{d\la^2}
=\int_{\pa B_1}\la^{1}\sum_{j=0}^2e_j\Delta_\theta^2\frac{d^2 u^\la}{d\la^2}\frac{d^2u^\la}{d\la^2}\\
&=\int_{\pa B_1}\la^{j+1}\Delta_\theta\frac{d^ju^\la}{d\la^j}\Delta_\theta\frac{d^2u^\la}{d\la^2}
\endaligned\ee
and
\be\nonumber\aligned
I_2(F_2):&=\int_{\pa B_1}\Big(7\Delta_\theta F_2(u^\la)-\la\frac{d}{d\la}\Delta_\theta^2 F_2(u^\la)\Big)\frac{du^\la}{d\la}\\
&=\int_{\pa B_1}\Big(7\Delta_\theta \big(\sum_{j=0}^2e_j\la^j\Delta^2_\theta\frac{d^ju^\la}{d\la^j}\big)-\la\frac{d}{d\la}\Delta_\theta^2 \big(\sum_{j=0}^2e_j\la^j\Delta^2_\theta\frac{d^ju^\la}{d\la^j}\big)\Big)\frac{du^\la}{d\la}\\
&=\sum_{j=1}^3\int_{\pa B_1}e_{0j}\Delta_\theta\frac{d^ju^\la}{d \la^j}\Delta_\theta\frac{du^\la}{d \la},
\endaligned\ee
where
\be\nonumber\aligned
e_{03}=-e_2,e_{02}=5e_2-e_1,e_{01}=6e_1-e_0
\endaligned\ee
and $e_j$ is defined in \eqref{12ej}.
\vskip0.1in
Recall the  formulas in the Lemma \ref{12fd1}-\ref{12fd2}, this time we regard that $f=\Delta_\theta u^\la$, then we can obtain that
\be\nonumber\aligned
I_1(F_2)+I_2(F_2)=&\int_{\pa B_1} C_1\la(\Delta_\theta\frac{du^\la}{d\la})^2+C_2\la^3(\Delta_\theta\frac{d^2u^\la}{d\la^2})^2\\
&\;\;+\frac{d}{d\la}\int_{\pa B_1}\Big(\sum_{0\leq i,j\leq 1}C_{i,j}\la^{i+j}\Delta_\theta\frac{d^iu^\la}{d\la^i}\Delta_\theta\frac{d^ju^\la}{d\la^j}\Big),
\endaligned\ee
where
\be\label{12c12}\aligned
C_1=-2e_0+6e_1-8e_2,C_2=2e_2
\endaligned\ee
and $C_{i,j}$ are determined by $e_j$,  we don't need  to outline the precise expressions of them  since the constant  will not  affect   our applying the monotonicity formula.

\subsection{The integral corresponding to  the operator $\Delta^3_\theta F_3(u^\la)$}

Recall the  formulas in the Lemma \ref{12fd1}-\ref{12fd2}
and  set $f=\nabla_\theta\Delta_\theta u^\la$,  then we can obtain that

\be\nonumber\aligned
I(F_3):&=\int_{\pa B_1}\Big(\la \Delta^3_\theta F_3(u^\la)\frac{d^2u^\la}{d\la^2}+7\Delta^3_\theta F_3(u^\la)\frac{du^\la}{d\la}-\la\frac{d \Delta^3_\theta F_3(u^\la)}{d\la}\frac{du^\la}{d\la}\Big)\\
&=\int_{\pa B_1}-\la\nabla_\theta\Delta_\theta u^\la \nabla_\theta\Delta_\theta\frac{d^2u^\la}{d\la^2}
-7\nabla_\theta\Delta_\theta u^\la\nabla_\theta\Delta_\theta \frac{du^\la}{d\la}
+\la|\nabla_\theta\Delta_\theta u^\la|^2\\
&=\frac{d}{d\la}\Big[\int_{\pa B_1}-\la \ds\Delta_{\theta} u^\la\ds\Delta_{\theta}\frac{du^\la}{d\la}
-3\big(\ds\Delta_{\theta}u^\la\big)^2\\
&\quad+2\la\int_{\pa B_1}\big( \ds\Delta_{\theta} \frac{du^\la}{d\la}\big)^2\Big].
\endaligned\ee

\label{12ccccc}\subsection{The monotonicity formula and the proof of Theorem \ref{12mono1} }

We sum up the terms  $\overline{E}_{d_1}$ and $\overline{E}_{d_2}$. Then
\be\nonumber\aligned
&E(\la,x,u):=
\int_{B_1}\frac{1}{2}|\Delta^2 u^\la|^2-\frac{1}{p+1}|u^\la|^{p+1}\\
&\quad+\int_{\pa B_1}\Big(\sum_{i,j\geq0,i+j\leq7}C_{i,j}^{0}\la^{i+j}\frac{d^iu^\la}{d\la^i}\frac{d^ju^\la}{d\la^j}
+\sum_{i,j\geq0,i+j\leq5}C_{i,j}^{1}\la^{i+j}\ds\frac{d^i u^\la}{d\la^i}\ds\frac{d^ju^\la}{d\la^j}\\
&\quad+\sum_{i,j\geq0,i+j\leq3}C_{i,j}^{2}\la^{i+j}\Delta_{\theta}\frac{d^iu^\la}{d\la^i}\Delta_{\q}\frac{d^ju^\la}{d\la^j}\\
&\quad+\sum_{i,j\geq0,i+j\leq1}C_{i,j}^{3}\la^{i+j}\ds\Delta_{\q}\frac{d^iu^\la}{d\la^i}\ds\Delta_{\q}\frac{d^ju^\la}{d\la^j}\Big),\\
\endaligned\ee
where the constant $C_{i,j}^k$ can be determined by the calculation of  $\overline{E}_{d_1}$ and $\overline{E}_{d_2}$ in the above three subsections.
Then we obtain the Theorem \ref{12mono1}.


\section{The desired monotonicity formula: The proof of Theorem \ref{12mono2}}

In this section, we construct the desired monotonicity formula  via the  blow-down analysis. We start from the Lemma \ref{12monoi} in the previous section. Firstly, we have
\bl If $p>\frac{n+8}{n-8}$, then
\be\nonumber
\sum_{l=1}^2(C_l+c_l)\la^{2l-1}(\Delta_\theta\frac{d^l u^\la}{d\la^l})^2\geq 0.
\ee
\el
\bp We known from \eqref{12c12}, \eqref{12ej} and \eqref{12abc} that
\be\nonumber\aligned
C_1=-6k^2+(-72+6n)k-178+30n;\;\;C_2+c_2=8,
\endaligned\ee
where $k=:\frac{8}{p-1}$. By this natation  we observe  that $p>\frac{n+8}{n-8}$ is  equivalent to $0<k<\frac{n-8}{2}$.
By finding  the roots (denoted  by  $r_1(n),r_2(n)$) of the equation $$-6k^2+(-72+6n)k-178+30n=0$$ about variable $k$, we get that
\be\nonumber\aligned
r_1(n):&=\frac{1}{2}n-6-\frac{1}{6}\sqrt{9n^2-36n+228}\\
&=\frac{1}{6}(3n-36-\sqrt{9n^2-36n+228})\\
&=\frac{1}{6}\frac{(3n-36)^2-(9n^2-36n+228)}{3n-36+\sqrt{9n^2-36n+228}}\\
&=\frac{-30n+178}{3n-36+\sqrt{9n^2-36n+228}}<0
\;\;\hbox{for}\;\;n\geq6
\endaligned\ee
and
\be\nonumber\aligned
r_2(n):=\frac{1}{2}n-6+\frac{1}{6}\sqrt{9n^2-36n+228}>\frac{1}{2}(n-8),
\endaligned\ee
therefore we obtain that $C_1>0$ if $0<k<\frac{n-8}{2}$. Recall that $c_2=2\alpha-2\beta-4>0$,
then  the conclusion follows.
\ep

\vskip0.2in

\bt\label{12BbBb} If $\frac{n+8}{n-8}<p<p_{c}(n)$, then there exist constants $b_{i,j}$ such that

\be\nonumber\aligned
\sum_{s=1}^3(B_s+b_s)\la^{2s-1}(\nabla_\theta\frac{d^su^\la}{d\la^s})^2
\geq\frac{d}{d\la}\Big(\sum_{0\leq i,j\leq2,i+j\leq3}b_{i,j}\nabla_\theta\frac{d^iu^\la}{d\la^i}\nabla_\theta\frac{d^ju^\la}{d\la^j}\Big).
\endaligned\ee

\et
\bp

To see this, from \eqref{12B123}, \eqref{12tj} and \eqref{12abc}, we get that
\be\nonumber\aligned
{Bb}_1:=B_1+b_1=&6k^4+(144-12n)k^3+(6n^2-204n+994)k^2\\
&\;\;+(60n^2-850n+2732)k+94n^2-1012n+2644,
\endaligned\ee
and
\be\nonumber\aligned
Bb_2:&=B_2+b_2=-38k^2+(-292+38n)k-6n^2+132n-544,\\
Bb_3:&=B_3+b_3=8.
\endaligned\ee
\bl\label{12Bb1} If $p>\frac{n+8}{n-8}$, then we have $Bb_1>0$.\el
\bp To show this, we first see that
if $n\in[9,29]$, under the condition $0<k<\frac{n-8}{2}$, by a direct  case by case calculation we can prove that
$Bb_1>0$.

For $n\geq29$, let us introduce the transform $k=\frac{n-8}{2}a$, hence $0<a<1$, then we have
\be\nonumber\aligned
{Bb}_1=&\frac{3}{8}(n-8)^4a^4-\frac{3}{2}(n-12)(n-8)^3a^3+\frac{1}{2}(3n^2-102n+497)(n-8)^2a^2\\
&\;\;+(n-8)(30n^2-425n+1366)a+94n^2-1012n+2644\\
=&f_4(n)a^4-f_3(n)a^3+f_2(n)a^2+f_1(n)a+f_0(n).
\endaligned\ee
We can see that if $n\geq29$, then $f_j(n)>0$ for $j=1,2,3,4$.
Since $0<a<1$, we have
\be\nonumber\aligned
{Bb}_1&=f_4(n)a^4-f_3(n)a^3+f_2(n)a^2+f_1(n)a+f_0(n)\\
&\geq f_4(n)a^4+(f_2(n)+f_1(n)+f_0(n)-f_2(n))a^3\\
&=f_4(n)a^4+(-1596+9n^3-\frac{261}{2}n^2+738n)a^3.\\
\endaligned\ee
By a direct calculation we can show that
\be\nonumber\aligned
-1596+9n^3-\frac{261}{2}n^2+738n>0\;\;\hbox{if}\;\;n\geq6.
\endaligned\ee
Hence we derive that $Bb_1>0$ when $n\geq29$. \ep

\bl\label{12Bb2} Assume that $\frac{n+8}{n-8}<p<p_{c}(n)$ and $n\geq9$, then $Bb_2>0$ except for $n=17,18$.
\el
\bp Firstly, we recall that $\frac{n+8}{n-8}<p<p_{c}(n)$ then we have $\frac{n-10}{2}-\sqrt{n}<k<\frac{n-8}{2}$ for $n\geq18$.
By solving   the equation  $Bb_2=0$, we  find the roots
\be\nonumber\aligned
r_1(n):=\frac{1}{2}n-\frac{73}{19}-\frac{1}{38}\sqrt{133n^2-532n+664},\\
r_2(n):=\frac{1}{2}n-\frac{73}{19}+\frac{1}{38}\sqrt{133n^2-532n+664}.
\endaligned\ee
Notice that $r_1(n)<\frac{n-10}{2}-\sqrt{n}$ is equivalent to
\be\nonumber\aligned
133n^2-1976n-3344\sqrt{n}-1292>0,
\endaligned\ee
the above inequality holds  whenever  $n\geq21$.
The strict inequality $r_2(n)>\frac{n-8}{2}$ can be proved  by a direct calculation.

Therefore, the conclusion holds when $n\geq21$.
For the remaining case  $9\leq n\leq20$ except for $n=17,18$, we can show  that the conclusion also holds.\ep

Combining with Lemmas \ref{12Bb1} and \ref{12Bb2}, we immediately get the following lemma.
\bl
For $n\in[9,+\infty),n\in  {\mathbb N}^+$, except for $n=17,18$, we have that
\be\nonumber\aligned
\sum_{s=1}^3(B_s+b_s)\la^{2s-1}(\nabla_\theta\frac{d^su^\la}{d\la^s})^2
\geq0.
\endaligned\ee
\el
Next, we have to consider the  remaining  case when $n=17,18$.
In view of Lemmas \ref{12fd1}-\ref{12fd2}, we have the following differential identity (denote that $f':=\nabla_\theta\frac{du^\la}{d\la}$):
\be\label{12biden}\aligned
\sum_{s=1}^3&(B_s+b_s)\la^{2s-1}(\nabla_\theta\frac{d^su^\la}{d\la^s})^2
=Bb_1\la(f')^2+(Bb_2+4Bb_3)\la^3(f'')^2\\
&\;\;+Bb_3\la(\la^2 f'''+2\la f'')^2+\frac{d}{d\la}\big(
-2Bb_3\la^4(f'')^2
\big)\\
&\geq(Bb_2+4Bb_3)\la^3(f'')^2+\frac{d}{d\la}\big(
-2Bb_3\la^4(f'')^2
\big).\\
\endaligned\ee
A direct calculation shows that under the condition $\frac{n+8}{n-8}<p<p_c(n)$ and $n=18$, we have that $Bb_2+4Bb_3>0$ (this way fails when $n=17$).
Hence, when $n=18$, we get that
\be\nonumber\aligned
\sum_{s=1}^3(B_s+b_s)\la^{2s-1}(\nabla_\theta\frac{d^su^\la}{d\la^s})^2
&\geq(Bb_2+4Bb_3)\la^3(f'')^2+\frac{d}{d\la}\big(-2Bb_3\la^4(f'')^2\big)\\
&\geq\frac{d}{d\la}\big(-2Bb_3\la^4(f'')^2\big).
\endaligned\ee
We need different method to handel  with $n=17$ in the following.
Notice that by the mean value inequality we have the following differential identity:
\be\nonumber\aligned
&\epsilon B_{10}\la(f')^2+B_{30}\la(\la^2f'''+2\la f'')^2
\geq -2\sqrt{\epsilon\cdot B_{30}\cdot B_{10}}(\la^3 f'''f'+2\la f'' f')\\
&=2\sqrt{\epsilon\cdot B_{30}\cdot B_{10}}(\la^3(f'')^2-\la(f')^2)+2\sqrt{\epsilon\cdot B_{30}\cdot B_{10}}
\frac{d}{d\la}\Big(\la^3 f''f'-\frac{1}{2}\la^2(f')^2\Big),
\endaligned\ee
here we have used the Lemmas \ref{12fd1}-\ref{12fd2};    $\epsilon\in (0,1)$ is to be determined later.

Combining with \eqref{12biden}, we have
\be\label{12n=18}\aligned
\sum_{s=1}^3&(B_s+b_s)\la^{2s-1}(\nabla_\theta\frac{d^su^\la}{d\la^s})^2
=Bb_1\la(f')^2+(Bb_2+4Bb_3)\la^3(f'')^2\\
&\;\;+Bb_3\la(\la^2 f'''+2\la f'')^2+\frac{d}{d\la}\big(-2Bb_3\la^4(f'')^2\big)\\
&\geq\Big(B_{20}+4B_{30}+2\sqrt{\epsilon\cdot B_{30}\cdot B_{10}}\Big)\la^3(f'')^2\\
&\;\;+\Big((1-\epsilon)B_{10}-2\sqrt{\epsilon\cdot B_{30}\cdot B_{10}}\Big)\la(f')^2\\
&\;\;+2\sqrt{\epsilon\cdot B_{30}\cdot B_{10}}
\frac{d}{d\la}\Big(\la^3 f''f'-\frac{1}{2}\la^2(f')^2\Big).
\endaligned\ee
To avoid confusion, we denote that $B_{j0}:=Bb_j$ for $j=1,2,3$ in the above and following.
Thus, our conclusion holds  once  we can prove
\be\label{12Bmean}
B_{20}+4B_{30}+2\sqrt{\epsilon\cdot B_{30}\cdot B_{10}}\geq0,
\ee
and
\be\label{12epsi}
(1-\epsilon)B_{10}-2\sqrt{\epsilon\cdot B_{30}\cdot B_{10}}\geq0.
\ee
To make sure \eqref{12epsi} holds, we select $\epsilon\in(0,1)$ satisfying
\be\nonumber
\frac{4\epsilon\cdot B_{30}}{(1-\epsilon)^2}\leq\min_{0\leq k\leq\frac{n-8}{2}} B_{10}.
\ee
Notice that when $n=17$, $\min_{0\leq k\leq\frac{n-8}{2}\mid_{n=17}} B_{10}=12606$, hence $\epsilon\leq0.9508$.
Thus we select that $\epsilon=0.9508$. Now we consider the inequality \eqref{12Bmean}.  Note that
 \be\nonumber \aligned
&\hbox{if}\;\;B_{20}+4B_{30}\geq0,\;\;\hbox{then the
inequality \eqref{12Bmean}    holds immediately;}\\
&\hbox{if}\;\;B_{20}+4B_{30}<0,\;\;\hbox{then \eqref{12Bmean} is  equivalent to}
\endaligned\ee
\be\label{12Bmean2}\;\;4\epsilon\cdot B_{30}\cdot B_{10}-(B_{20}+B_{30})^2\geq0.\ee
By a direct calculation we can prove that when $n=17$, under the condition $\frac{n+8}{n-8}<p<p_c(n)$, i.e., $\max\{0,R_1(n)\}<k<\frac{n-8}{2}$ ($R_1(n)$ see \eqref{LWZ-666}), we have that
\be\label{12n=18}\aligned
&\sum_{s=1}^3(B_s+b_s)\la^{2s-1}(\nabla_\theta\frac{d^su^\la}{d\la^s})^2\geq\Big(B_{20}+4B_{30}+2\sqrt{\epsilon\cdot B_{30}\cdot B_{10}}\Big)\la^3(f'')^2\\
&+\Big((1-\epsilon)B_{10}-2\sqrt{\epsilon\cdot B_{30}\cdot B_{10}}\Big)\la(f')^2\\
&+2\sqrt{\epsilon\cdot B_{30}\cdot B_{10}}
\frac{d}{d\la}\Big(\la^3 f''f'-\frac{1}{2}\la^2(f')^2\Big)\\
&\geq2\sqrt{\epsilon\cdot B_{30}\cdot B_{10}}
\frac{d}{d\la}\Big(\la^3 f''f'-\frac{1}{2}\la^2(f')^2\Big).\\
\endaligned\ee
Summing up,  for all cases, we have proved  Theorem \ref{12BbBb}.
\ep
\vskip0.2in

\bt\label{12AaAa} If $\frac{n+8}{n-8}<p<p_c(n)$ and $n\geq9$, then there exist constants $a_{i,j}$ such that
\be\nonumber\aligned
\sum_{j=1}^4(A_j+a_j)\la^{2j-1}(\frac{d^ju^\la}{d\la^j})^2\geq
\frac{d}{d\la}\Big(\sum_{0\leq i,j\leq2,i+j\leq5}a_{i,j}\frac{d^iu^\la}{d\la^i}\frac{d^ju^\la}{d\la^j}\Big).
\endaligned\ee
\et
  Firstly, from \eqref{12d1d2d3}, \eqref{12abc} and \eqref{12kj} we have that
\be\nonumber\aligned
&Aa_1:=A_1+a_1=-4k^6+(-88+12n)k^5+(-12n^2+208n-860)k^4\\
&\;\;+(4n^3-152n^2+1544n-4304)k^3
+(32n^3-776n^2+5408n-11196)k^2\\
&\;\;+(92n^3-1576n^2+8428n-14040)k+64n^3-940n^2+4368n-6372,
\endaligned\ee

\be\label{12Aa2}\aligned
&Aa_2:=A_2+a_2=28k^4+(464-56n)k^3+(32n^2-656n+2668)k^2\\
&+(-4n^3+232n^2-2364n+6456)k-16n^3+380n^2-2640n+5556,
\endaligned\ee

\be\nonumber\aligned
Aa_3:=A_3+a_3=-28k^2+(-216+28n)k-4n^2+96n-408
\endaligned\ee
and $Aa_4:=A_4+a_4=4$.

We separate the proofs into several Lemmas.

\bl\label{12A1a1} If $p>\frac{n+8}{n-8}$ and $n\geq9$, then $A_1+a_1>0$.
\el
\bp To see this, in fact from \eqref{12d1d2d3}, \eqref{12abc} and \eqref{12kj} we have that
\be\nonumber\aligned
A_1=&-2k^6+(-72+6n)k^5+(-6n^2+174n-818)k^4+(2n^3-132n^2+\\
&1492n-4272)k^3+(30n^3-768n^2+5412n-11222)k^2+(94n^3-1596n^2\\
&+8486n-14088)k+66n^3-954n^2+4398n-6390,
\endaligned\ee
and
\be\nonumber\aligned
a_1=&-2k^6+(-16+6n)k^5+(-6n^2+34n-42)k^4+(2n^3-20n^2+52n\\
&-32)k^3+(2n^3-8n^2-4n+26)k^2+(-2n^3+20n^2-58n+48)k\\
&-2n^3+14n^2-30n+18.
\endaligned\ee
It is can be observed that $-3,-1$ and $n-3,n-5$
are the roots of $A_1=0$, hence we have that
\be\nonumber\aligned
A_1=(k+3)(k+1)\big(k-(n-3)\big)\big(k-(n-5)\big)(-2k^2+(2n-48)k+22n-142),
\endaligned\ee
from this, it is not difficult to see that
$A_1>0$ if $0<k<\frac{n-8}{2}$.
For $a_1=0$, we have observed that $-1,-1,1,n-3,n-3,n-1$ are the roots, hence
\be\nonumber\aligned
a_1=(k+1)^2(k-1)\big(k-(n-3)\big)^2\big(k-(n-1)\big),
\endaligned\ee
thus $a_1>0$ if $1<k<n-1$.
For $n\geq18$, we have that $\frac{n-10}{2}-\sqrt{n}>1$, then $a_1>0$ if $\frac{n+8}{n-8}<p$, i.e., $0<k<\frac{n-8}{2}$.
Once $a_1>0$, then $A_1+a_1>0$. For the case $n\in [9,17]$, we calculate case by case showing  that $A_1+a_1>0$ if $0<k<\frac{n-8}{2}$.
This completes  the proof of Lemma \ref{12A1a1}. \ep

\vskip0.12in

\bl\label{12A3a3} If $\frac{n+8}{n-8}<p<p_c(n)$ and $n\geq9$, then $A_3+a_3>0$.
\el
\bp  By solving  the roots of $A_3+a_3=0$, we get that
\be\nonumber\aligned
r_{10}(n):=\frac{1}{2}n-\frac{27}{7}-\frac{1}{14}\sqrt{21n^2-84n+60},\\
r_{20}(n):=\frac{1}{2}n-\frac{27}{7}+\frac{1}{14}\sqrt{21n^2-84n+60}.
\endaligned\ee
Notice that $r_{10}(n)<\frac{n-10}{2}-\sqrt{n}$ is  equivalent to $$21n^2-85n-32\sqrt{n}-196>0.$$
The above inequality holds if $n\geq7$. The root $r_{20}(n)>\frac{n-8}{2}$ can be seen immediately.
Thus, for $n\geq18$ we know that  $A_3+a_3>0$ if $R_1(n)<k<\frac{n-8}{2}$.
However, when $n\in[9,17]$, we may compute case by case to show that $A_3+a_3>0$  if  $\frac{n+8}{n-8}<p<p_c(n)$. This finishes  the proof of Lemma \ref{12A3a3}.
\ep

\bl\label{12A2a2}
If $\frac{n+8}{n-8}<p<p_c(n)$, for $n\in[9,13]$ or $n\geq21$, we have that $A_2+a_2>0$.
\el

\bp To prove Lemma \ref{12A2a2}, we introduce the transformation  $n:=t^2, k:=\frac{n-10}{2}-a\cdot t$,
hence $0<a<1$.
From \eqref{12Aa2}, we have that
\be\nonumber\aligned
A_2+a_2&=-524+\frac{3}{4}t^8+(-10a^2-6)t^6-24at^5+(28a^4+40a^2-5)t^4\\
&\;\;+(96a^3+96a)t^3+(-92a^2+68)t^2-576at\\
&\geq-524+\frac{3}{4}t^8-16t^6-24t^5-5t^4-24t^2-576t>0\;\hbox{if}\;n=t^2\geq28.
\endaligned\ee
For the cases  $n\in[9,13]$ and $n\in[21,27]$, we may calculate  case by case to  show that
$A_2+a_2>0$ if $\frac{n+8}{n-8}<p<p_c(n)$. This finishes  the proof of Lemma \ref{12A2a2}.\ep
\vskip0.1in

Combining with Lemmas \ref{12A1a1}, \ref{12A3a3} and \ref{12A2a2}, we have that
\bl If $\frac{n+8}{n-8}<p<p_c(n)$, for $n\in[9,13]$ and $n\geq21$, we have
\be\nonumber\aligned
\sum_{j=1}^4(A_j+a_j)\la^{2j-1}(\frac{d^ju^\la}{d\la^j})^2\geq0.
\endaligned\ee
\el

\br
To establish the Theorem \ref{12AaAa}, we need to deal   with the case $n\in[14,20]$ by some other methods.
\er

Now we turn to the case $n\in[14,20]$.
To deal with  these cases, firstly  we establish  the following differential identity (denote that $f'=\frac{du^\la}{d\la}$):
\be\label{12d1234}\aligned
&\sum_{j=1}^4(A_j+a_j)\la^{2j-1}(\frac{d^ju^\la}{d\la^j})^2:=\sum_{j=1}^4d_j\la^{2j-1}(\frac{d^ju^\la}{d\la^j})^2\\
&=d_1\la(f')^2+(d_2+4d_3)\la^3(f'')^2+d_3\la(\la^2f'''+2\la f'')^2\\
&\;\;+d_4\la^7(f'''')^2+\frac{d}{d\la}(-2d_3\la^4(f'')^2),
\endaligned\ee
where we have used Lemma \ref{12fd2}.
The term $d_2+4d_3$   is not   nonnegative. Now we invoke the mean value inequality.
For the parameter $x\in[0,1]$ to be determined later, we have that
\be\label{12d14}\aligned
& d_1\la(f')^2+d_4\la^7(f'''')^2=x\cdot d_1\la(f')^2+d_4\la^7(f'''')^2+(1-x)\cdot d_1\la(f')^2\\
&\geq2\sqrt{x\cdot d_1d_4}\la^4f''''f'+(1-x)d_1\la(f')^2\\
&=2\sqrt{x\cdot d_1d_4}\Big(-12\la(f')^2+6\la^3(f'')^2\Big)+(1-x)d_1\la(f')^2\\
&\;\;+\frac{d}{d\la}\Big(2\sqrt{x\cdot d_1d_4}\big(\la^4 f'''f'-\frac{1}{2}\la^4(f'')^2-4\la^3f''f'+6\la^2(f')^2\big)\Big)\\
&=\Big((1-x)d_1-24\sqrt{x\cdot d_1\cdot d_4}\Big)\la(f')^2+12\sqrt{x\cdot d_1\cdot d_4}\la^3(f'')^2\\
&\;\;+\frac{d}{d\la}\Big(\sum_{1\leq i,j\leq3,i+j\leq4}d_{i,j}\la^{i+j}f^{(i)}f^{(j)}\Big),
\endaligned\ee
where we have used Lemma \ref{12fd2}.
The constants in the derivative terms,   namely, $d_{i,j}$ can be determined but we do not need  the exactly expressions. In particular,  $d_{i,j}$ may be changed in the following derivation, but we still denote as $d_{i,j}$.
Combine \eqref{12d1234} and \eqref{12d14}, we  obtain that
\be\label{12dnnn}\aligned
&\sum_{j=1}^4(A_j+a_j)\la^{2j-1}(\frac{d^ju^\la}{d\la^j})^2:=\sum_{j=1}^4d_j\la^{2j-1}(\frac{d^ju^\la}{d\la^j})^2\\
&=d_1\la(f')^2+(d_2+4d_3)\la^3(f'')^2+d_3\la(\la^2f'''+2\la f'')^2+d_4\la^7(f'''')^2\\
&\;\;+\frac{d}{d\la}(-2d_3\la^4(f'')^2),\\
&\geq\Big((1-x)d_1-24\sqrt{x\cdot d_1\cdot d_4}\Big)\la(f')^2
+(d_2+4d_3+12\sqrt{x\cdot d_1\cdot d_4})\la^3(f'')^2\\
&\;\;+d_3\la(\la^2f'''+2\la f'')^2
+\frac{d}{d\la}\Big(\sum_{1\leq i,j\leq3,i+j\leq4}d_{i,j}\la^{i+j}f^{(i)}f^{(j)}\Big)\\
&\geq\Big((1-x)d_1-24\sqrt{x\cdot d_1\cdot d_4}\Big)\la(f')^2
+(d_2+4d_3+12\sqrt{x\cdot d_1\cdot d_4})\la^3(f'')^2\\
&\;\;+\frac{d}{d\la}\Big(\sum_{1\leq i,j\leq3,i+j\leq4}d_{i,j}\la^{i+j}f^{(i)}f^{(j)}\Big).\\
\endaligned\ee
If
\be\label{12dx}\aligned
(1-x)d_1-24\sqrt{x\cdot d_1\cdot d_4}>0
\endaligned\ee
and
\be\label{12dd}\aligned
d_2+4d_3+12\sqrt{x\cdot d_1\cdot d_4}\geq0
\endaligned\ee
hold simultaneously, then we have Theorem \ref{12AaAa}.
By this method, we have the following

\bl\label{12n=1420} If $\frac{n+8}{n-8}<p<p_c(n)$, and $n\in[14,20]$ except for $n=17$, then there exist constant $d_{i,j}$ such that
\be\nonumber\aligned
\sum_{j=1}^4(A_j+a_j)\la^{2j-1}(\frac{d^ju^\la}{d\la^j})^2
\geq\frac{d}{d\la}\Big(\sum_{1\leq i,j\leq3,i+j\leq4}d_{i,j}\la^{i+j}\frac{d^iu^\la}{d\la^i}\frac{d^ju^\la}{d\la^j}\Big).
\endaligned\ee
\el
\bp  To prove this  Lemma, we just need to check that the conditions \eqref{12dx} and \eqref{12dd} hold simultaneously  when  $n\in[14,20]$ except for $n=17$.  We perform these case by case.
To make sure that \eqref{12dx} hold, we need to  select that $x\in[0,1]$ satisfying
\be\label{12dxx}
\frac{24^2d_4\cdot x}{(1-x)^2}<\min_{0\leq k\leq\frac{n-8}{2}}d_1.
\ee
On the other hand, condition \eqref{12dd} can be obtained by the following two cases:
\be\label{12ddd}\aligned
&\hbox{If}\;\;d_2+4d_3\geq0,\;\;\hbox{then}\;\;\eqref{12dd} \;\;\hbox{holds immediately};\\
&\hbox{If}\;\;d_2+4d_3<0,\;\;\hbox{then}\;\;\eqref{12dd} \;\;\hbox{is equivalent to}\;\;12^2x\cdot d_1\cdot d_4-(d_2+4d_3)^2>0.
\endaligned\ee
For the  simplicity, we denote that $d:=d(k,n,x)=12^2x\cdot d_1\cdot d_4-(d_2+4d_3)^2$.
Now we turn to consider the inequalities \eqref{12dxx} and \eqref{12ddd}. Let
\be\label{LWZ-222} R_1(n)=\frac{n-10}{2}-d(n),\ee
where $d(n)$ is defined in \eqref{LWZ=333}. Note  that $p<p_c(n)$,  then $\max\{R_1(n),0\}<k$. Recall that $k:=\frac{8}{p-1}$.

\vskip0.1in
For $n=14$, $\min_{0\leq k\leq\frac{n-8}{2}\mid_{n=14}}d_1=46156$, then from \eqref{12dxx} we get that $x\leq0.8001464380$, thus we select
$x=0.8001464380$.  Thus,  $d_2+4d_3<0$ if $0<k<0.02572910109$ and $d(k,n,x)_{n=14,x=0.8001464380}>0$ if $0<k< 0.7919464848$. Hence,  by selecting  the parameter
$x=0.8001464380$, we can make sure that inequalities \eqref{12dxx} and \eqref{12ddd} hold simultaneously.  Hence,  Lemma \ref{12n=1420} holds  for the case $n=14$.

\vskip0.1in
Similarly we may prove  Lemma \ref{12n=1420} for $n=15, 16, 17, 18, 19$. We omit these details. But we would like to give the proof for the case of $n=20$ for reader's convenience.

\vskip0.1in
Let $n=20$. Then $\min_{0\leq k\leq\frac{n-8}{2}\mid_{n=20}}d_1=216988$. From \eqref{12dxx} we get that $x\leq0.9021282144$. Thus we may select
$x=0.9021282144$.  Hence, $d_2+4d_3<0$ if $.9244642513\approx R_1(n=20)<k<1.026523007$ and $d(k,n,x)_{n=20,x=0.9021282144}>0$
 if $.9244642513\approx R_1(n=20)<k< 1.894875455$. So, by selecting  the parameter
$x=0.9021282144$, then  both  inequalities \eqref{12dxx} and \eqref{12ddd} hold simultaneously. Therefore,    Lemma \ref{12n=1420} holds when $n=20$.
\ep


\vskip0.23in

\br\label{12re17} For the case $n=17$, this method of Lemma \ref{12n=1420}  does not work. But if $k>0.02175341614$, then we have Lemma \ref{12n=1420} hold.
Since, when   $n=17$, $$\min_{0\leq k\leq\frac{n-8}{2}\mid_{n=17}}d_1=110656,$$ then from \eqref{12dxx} we get that $x\leq0.8657397553$, thus we select
$x=0.8657397553$. It follows that $d_2+4d_3<0$ if $0<k<0.5256119817$ and that $$d(k,n,x)_{n=17,x=0.8657397553}>0  \hbox{  if } 0.02175341614<k< 1.358050900.$$ Thus, by selecting the parameter
$x=0.8657397553$ and  $k>0.02175341614$, then the  inequalities \eqref{12dxx} and \eqref{12ddd} hold simultaneously, hence the  Lemma \ref{12n=1420} holds.
\er
\br We note that the dimension $n=17$ is critical dimension when we define the Joseph-Lundgren exponent in \eqref{12pcn}, it is the most thorny  dimension when we establish
Theorem \ref{12AaAa}.
\er
But by getting more information from \eqref{12dnnn}, we still have the following

\bl\label{12n=17} If $\frac{n+8}{n-8}<p<p_c(n)$  and $n=17$, then there exist constant $d_{i,j}$ such that
\be\nonumber\aligned
\sum_{j=1}^4(A_j+a_j)\la^{2j-1}(\frac{d^ju^\la}{d\la^j})^2
\geq\frac{d}{d\la}\Big(\sum_{1\leq i,j\leq3,i+j\leq4}d_{i,j}\la^{i+j}\frac{d^iu^\la}{d\la^i}\frac{d^ju^\la}{d\la^j}\Big).
\endaligned\ee
\el
\bp In view of Remark \ref{12re17}, we only need to consider the case  of $0<k<0.04$.
From \eqref{12dnnn}, we discard the the term $d_3\la(\la^2f'''+2\la f'')^2$ (which is nonnegative term hence a ''good'' term) directly in our estimate, now we ''pick up'' and make full use of this term. To achieve this, let us select parameters $x_1,x_2,y\in[0,1]$ whose exact values are  to be determined later. We have
\be\nonumber\aligned
&y\cdot d_1\cdot\la(f')^2+d_3\cdot\la(\la^2f''+2\la f')^2\\
&=x_1\cdot y\cdot d_1\cdot\la(f')^2+d_3\cdot\la(\la^2f''+2\la f')^2+(1-x_1)\cdot y\cdot d_1\cdot\la(f')^2\\
&\geq2\sqrt{x_1\cdot y\cdot d_1\cdot d_3}\la^3(f'')^2+\Big((1-x_1)\cdot y\cdot d_1-2\sqrt{x_1\cdot y\cdot d_1\cdot d_3}\Big)\la(f')^2\\
&\;\;+\frac{d}{d\la}\Big(2\sqrt{x_1\cdot y\cdot d_1\cdot d_3}(\la^3 f''f'-\frac{1}{2}\la^2(f')^2)\Big) \\
\endaligned\ee
and

\be\nonumber\aligned
& (1-y)\cdot d_1\la(f')^2+d_4\la^7(f'''')^2\\
&=x_2\cdot(1-y)\cdot d_1\la(f')^2+d_4\la^7(f'''')^2+(1-x_2)\cdot(1-y)\cdot d_1\la(f')^2\\
&\geq2\sqrt{x_2\cdot(1-y)\cdot d_1\cdot d_4}\la^4f''''f'+(1-x_2)(1-y)d_1\la(f')^2\\
&=2\sqrt{x_2(1-y)d_1d_4}\Big(-12\la(f')^2+6\la^3(f'')^2\Big)+(1-x_2)(1-y)d_1\la(f')^2\\
&\quad+\frac{d}{d\la}\Big(2\sqrt{x\cdot d_1d_4}\big(\la^4 f'''f'-\frac{1}{2}\la^4(f'')^2-4\la^3f''f'+6\la^2(f')^2\big)\Big)\\
&=\Big((1-x_2)(1-y)d_1-24\sqrt{x_2(1-y)\cdot d_1\cdot d_4}\Big)\la(f')^2\\
&\quad+12\sqrt{x_2(1-y)\cdot d_1\cdot d_4}\la^3(f'')^2+\frac{d}{d\la}\Big(\sum_{1\leq i,j\leq3,i+j\leq4}d_{i,j}\la^{i+j}f^{(i)}f^{(j)}\Big).
\endaligned\ee
Therefore from \eqref{12d1234}, combine with the  above two inequalities, we get that

\be\nonumber\aligned
&\sum_{j=1}^4(A_j+a_j)\la^{2j-1}(\frac{d^ju^\la}{d\la^j})^2:=\sum_{j=1}^4d_j\la^{2j-1}(\frac{d^ju^\la}{d\la^j})^2\\
&=d_1\la(f')^2+(d_2+4d_3)\la^3(f'')^2+d_3\la(\la^2f'''+2\la f'')^2\\
&\quad+d_4\la^7(f'''')^2+\frac{d}{d\la}(-2d_3\la^4(f'')^2)\\
&\geq(d_2+4d_3+2\sqrt{x_1\cdot y\cdot d_1\cdot d_3}+12\sqrt{x_2(1-y)\cdot d_1\cdot d_4})\la^3(f'')^2\\
&\quad+\Big((1-x_1)\cdot y\cdot d_1+(1-x_2)(1-y)d_1-2\sqrt{x_1\cdot y\cdot d_1\cdot d_3}\\
&\quad-24\sqrt{x_2(1-y)\cdot d_1\cdot d_4}\Big)\la(f')^2+\frac{d}{d\la}\Big(\sum_{1\leq i,j\leq3,i+j\leq4}d_{i,j}\la^{i+j}f^{(i)}f^{(j)}\Big).
\endaligned\ee
We may establish the Lemma \eqref{12n=17} if we have

\be\label{12d10703}\aligned
(1-x_1) y d_1+(1-x_2)(1-y)d_1-2\sqrt{x_1 y d_1 d_3}-24\sqrt{x_2(1-y) d_1 d_4}>0
\endaligned\ee
and
\be\label{12d20703}\aligned
d_2+4d_3+2\sqrt{x_1\cdot y\cdot d_1\cdot d_3}+12\sqrt{x_2(1-y)\cdot d_1\cdot d_4}\geq0
\endaligned\ee
hold simultaneously.
We will select proper parameters $x_1,x_2,y\in[0,1]$ to make sure the  inequalities \eqref{12d10703} and \eqref{12d20703} hold under the condition
$0<k<0.04$.
For simplicity, we denote that
\be\nonumber\aligned
f_1:&=\Big((1-x_1) y +(1-x_2)(1-y)\Big)^2d_1-4x_1yd_3-24^2x_2(1-y)d_4,\\
f_2:&=f_1^2-(4\times24)^2x_1yx_2(1-y)d_3d_4,\\
h_1:&=(d_2+4d_3)^2-4x_1yd_1d_3-12^2x_2(1-y)d_1d_4,\\
h_2:&=48^2x_1yx_2(1-y)d_3d_4d_1^2-h_1^2.\\
\endaligned\ee
Hence, \eqref{12d10703} is  equivalent to $f_1>0$ and $f_2>0$. As for  \eqref{12d20703} we observe that
\be\label{120704h}\aligned
&\hbox{if}\;\;h_1\leq0,\;\;\hbox{then}\;\; \eqref{12d20703} \;\;\hbox{holds immediately},\\
&\hbox{if}\;\;h_1>0,\;\;\hbox{then}\;\; h_2>0.
\endaligned\ee
How can we select the proper parameter $x_1,x_2,y$ ? Numerically, by considering  the end-point case, namely $y=0$
and $y=1$, then as in \eqref{12dxx} to determine $x_1,x_2$, we find that we may select   $y$ closing  to $0$
and $x_1,x_2$ closing  to $1$ . In fact, let us select that
$y=0.1,x_1=x_2=0.8$, we have that
\be\nonumber\aligned
h_1&=748.16k^8-25955.84k^7+2.5965440\cdot10^5k^6+1.8144000\cdot10^5k^5\quad\quad\\
&\quad-1.370264514\cdot10^7k^4+2.783143502\cdot10^7k^3+1.905355534\cdot10^8k^2\\
&\quad-2.946605150\cdot10^8k+1.316646282\cdot10^7,
\endaligned\ee
\be\nonumber\aligned
h_2&=2.689086\cdot10^{14}+3.883824251\cdot10^7k^{15}-5.597433856\cdot10^5k^{16}\\
&\;\;+1.045307292\cdot10^{16}k-1.881709367\cdot10^{16}k^5
+8.22093182\cdot10^{12}k^{10}\\
&\;\;+8.52692786\cdot10^{14}k^7+7.5489982\cdot10^{11}k^9-3.79196866\cdot10^{10}k^{12}\\
&\;\;+1.322360547\cdot10^{10}k^{13}
+1.138472739\cdot10^{17}k^3-1.062469518\cdot10^9k^{14}\\
&\;\;-8.42043516\cdot10^{11}k^{11}-1.943781276\cdot10^{16}k^4-3.114952088\cdot10^{14}k^8\\
&\;\;-8.750277873\cdot10^{16}k^2+4.523393231\cdot10^{15}k^6,
\endaligned\ee

\be\nonumber\aligned
f_1=&-0.1600k^6+4.6400k^5-31.6800k^4-93.2800k^3+556.6400k^2\quad\quad\quad\quad\quad\\
&\;\;+4947.5200k+2745.6000,
\endaligned\ee
\be\nonumber\aligned
f_2=&-40.20787k^8+31.6672860k^{10}+7.3936205\cdot10^6+9492.1751k^7\\
&\;\;+2.66134182\cdot 10^7k
+0.02560000k^{12}-1.48480312k^{11}-3.9194703\cdot10^5k^5\\
&\;\;+4.99650206\cdot10^6k^3
-7.8719259\cdot10^5k^4+2.75922586\cdot10^7k^2\\
&\;\;+18460.169k^6-264.138939k^9.
\endaligned\ee
We can verify that
$f_1>0$ if $0<k<13.82353260$; $f_2>0$ if $0<k<9.306459393$;
$h_1>0$ if $0<k<0.04606463340$; $h_2>0$ if $0<k<0.1757218049$.

Thus, when $0<k<0.04$, we have that $f_1,f_2,h_1,h_2$ all are positive. Hence,  combine with \eqref{120704h}, we
know that \eqref{12d10703} and \eqref{12d20703} hold. This  proves the Lemma \ref{12n=17}.
\ep


\section{Homogeneous stable solutions must be zero solution}

Firstly, we establish the representations of the  harmonic, biharmonic, triharmonic and quadharmonic  operators in terms of the  spherical coordinates.
We continue from  the Section 4. By a direct calculation we have
\be\nonumber\aligned
\Delta^4u=&(\pa_{rr}+\frac{n-1}{r}\pa_r)^4+\Delta_\theta^4(\pa_{rr}+\frac{n-1}{r}\pa_r)(r^{-8}u)
+\Delta_\theta\Big(
(\pa_{rr}+\frac{n-1}{r}\pa_r)^2\\
&\big(r^{-2}(\pa_{rr}+\frac{n-1}{r}\pa_r)u\big)+(\pa_{rr}+\frac{n-1}{r}\pa_r)^3(r^{-2}u)
+(\pa_{rr}+\frac{n-1}{r}\pa_r)\cdot\\
&\big(r^{-2}(\pa_{rr}+\frac{n-1}{r}\pa_r)^2u\big)+r^{-2}(\pa_{rr}+\frac{n-1}{r}\pa_r)^3u
\Big)\\
&+\Delta_\theta^2\Big(
(\pa_{rr}+\frac{n-1}{r}\pa_r)^2(r^{-4}u)+(\pa_{rr}+\frac{n-1}{r}\pa_r)\big(r^{-4}(\pa_{rr}+\frac{n-1}{r}\pa_r)u\big)
\Big)\\
&+(\pa_{rr}+\frac{n-1}{r}\pa_r)\big(r^{-2}(\pa_{rr}+\frac{n-1}{r}\pa_r)r^{-2}u\big)
+r^{-2}(\pa_{rr}+\frac{n-1}{r}\pa_r)\cdot\\
&\big(r^{-2}(\pa_{rr}+\frac{n-1}{r}\pa_r)u\big)
+r^{-2}(\pa_{rr}+\frac{n-1}{r}\pa_r)^2(r^{-2}u)+r^{-4}(\pa_{rr}+\frac{n-1}{r}\pa_r)^2u\\
&+\Delta_\theta^3\Big(
(\pa_{rr}+\frac{n-1}{r}\pa_r)(r^{-6}u)+r^{-2}(\pa_{rr}+\frac{n-1}{r}\pa_r)(r^{-4}u)\\
&+r^{-6}(\pa_{rr}+\frac{n-1}{r}\pa_r)u+r^{-4}(\pa_{rr}+\frac{n-1}{r}\pa_r)(r^{-2}u)
\Big).
\endaligned\ee
Then, let $u=r^{-k}w(\theta)$, recall that $k=\frac{8}{p-1}$, by a direct calculation, the function  $w$ satisfy
\be\label{12ww}\aligned
\Delta_\theta^4w-J_3\Delta_\theta^3w+J_2\Delta_\theta^2w-J_1\Delta_\theta w+J_0w=|w|^{p-1}w,
\endaligned\ee
where
\be\nonumber\aligned
J_0=k(k+2-n)(k+2)(k+4-n)(k+4)(k+6-n)(k+6)(k+8-n),
\endaligned\ee
\be\label{12J1}\aligned
-J_1&=k(k+2-n)(k+4)(k+6-n)(k+6)(k+8-n)\quad\quad\quad\quad\quad\quad\quad\quad \\
&+(k+2)(k+4-n)(k+4)(k+6-n)(k+6)(k+8-n)\\
&+k(k+2-n)(k+2)(k+4-n)(k+6)(k+8-n)\\
&+k(k+2-n)(k+2)(k+4-n)(k+4)(k+6-n),\\
\endaligned\ee
\be\label{12J2}\aligned
J_2=&(k+4)(k+6-n)(k+6)(k+8-n)+k(k+2-n)(k+6)(k+8-n)\\
+&(k+2)(k+4-n)(k+6)(k+8-n)+k(k+2-n)(k+4)(k+6-n)\\
+&(k+2)(k+4-n)(k+4)(k+6-n)+k(k+2-n)(k+2)(k+4-n),
\endaligned\ee
and
\be\label{12J3}\aligned
-J_3&=(k+6)(k+8-n)+(k+4)(k+6-n)+k(k+2-n)\\
&+(k+2)(k+4-n).
\endaligned\ee
Since $w\in W^{4,2}(\mathbb{S}^{n-1})\cap L^{p+1}(\mathbb{S}^{n-1})$, $r^{-\frac{n-8}{2}}w(\theta)\eta_\varepsilon(r)$
can be approximated  by $C_0^\infty(B_{4/\varepsilon}\setminus B_{\varepsilon/4})$ functions in
$W^{4,2}(B_{2/\varepsilon}\setminus B_{\varepsilon/2})\cap L^{p+1}(B_{2/\varepsilon}\setminus B_{\varepsilon/2})$.
Hence, here we may insert the    the stability condition of  $u$ and  choose a test function of the form $\varphi=r^{-\frac{n-8}{2}}w(\theta)\eta_\varepsilon(r)$.
Note that
\be\nonumber\aligned
\Delta^2\varphi&=r^{-\frac{n}{2}}\eta_\varepsilon(r)
\Big(
q(q+2-n)(q+2)(q+4-n)w(\theta)+\big((q+2)(q+4-n)\\
&+q(q+2-n)\big)\Delta_\theta w+\Delta^2_\theta w
\Big)+\sum_{j=1}^4c_jr^{-\frac{n}{2}+j}\eta_\varepsilon^{(j)}(r)w(\theta),
\endaligned\ee
where $q=\frac{n-8}{2}$, $c_j$ are some constants, $\eta^{(j)}:=\frac{d^j}{dr^j}\eta_\varepsilon$ for $j\geq1$ and $\eta_\varepsilon^{(0)}:=\eta_\varepsilon$.
Hence we have that
\be\label{12homo1}\aligned
&\int_{\R^n}|\Delta^2 \varphi|^2=\Big(\int_{0}^\infty r^{-1}\eta_\varepsilon^2(r)dr\Big)
\Big(\int_{\mathbb{S}^{n-1}}q_0w^2(\theta)+q_1|\nabla_\theta w|^2+q_2|\Delta_\theta w|^2\\
&+q_3|\nabla_\theta\Delta_\theta w|^2+|\Delta^2_\theta w|^2\Big)
+\int_{0}^\infty\sum_{1\leq i+j\leq8,i,j\geq0}c_{i,j}r^{i+j-1}\eta_\varepsilon^{(i)}\eta_\varepsilon^{(j)}dr\\
&\cdot\Big(
\int_{\mathbb{S}^{n-1}}w^2(\theta)+|\nabla_\theta w|^2+|\Delta_\theta w|^2
+|\nabla_\theta\Delta_\theta w|^2+|\Delta^2_\theta w|^2
\Big),
\endaligned\ee
where
\be\label{12q}\aligned
q_0&=\Big(q(q+2-n)(q+2)(q+4-n)\Big)^2,\\
q_1&=2\Big((q+2)(n-q-4)+q(n-q-2)\Big)q(q+2-n)(q+2)(q+4-n)\\
q_2&=\Big((q+2)(q+4-n)+q(q+2-n)\Big)^2+2q(q+2-n)(q+2)(q+4-n)\\
q_3&=2\Big((q+2)(n-q-4)+q(n-q-2)\Big)\\
\endaligned\ee
Substituting this into the stability condition for $u$, we get that
\be\label{12homo2}\aligned
p\big(\int_{\mathbb{S}^{n-1}}|w|^{p+1}d\theta\big)\cdot\big(\int_{0}^\infty r^{-1}\eta_\varepsilon^2(r)dr\big)
\leq\int_{\R^n}|\Delta^2 \varphi|^2.
\endaligned\ee
Notice that
\be\nonumber\aligned
&\int_0^\infty r^{-1}\eta_\varepsilon^2(r)dr\geq |\log \varepsilon|\rightarrow+\infty,\;\;\hbox{as}\;\; \varepsilon\rightarrow0,\\
&\int_0^\infty \sum_{1\leq i+j\leq8,i,j\geq0}r^{i+j-1}|\eta_\varepsilon^{(i)}\eta_\varepsilon^{(j)}|dr\leq C,\;\;\hbox{for any}\;\;i,j.
\endaligned\ee
Combine with \eqref{12homo1} and \eqref{12homo2}, we obtain that
\be\nonumber\aligned
p\int_{\mathbb{S}^{n-1}}|w|^{p+1}d\theta&\leq\int_{\mathbb{S}^{n-1}}q_0w^2(\theta)+q_1|\nabla_\theta w|^2+q_2|\Delta_\theta w|^2\\
&+q_3|\nabla_\theta\Delta_\theta w|^2+|\Delta^2_\theta w|^2,
\endaligned\ee
Combining this with \eqref{12ww}, we have the following estimate,
\be\label{12JPQ}\aligned
&\int_{\mathbb{S}^{n-1}}(p J_0-q_0)w^2(\theta)+(p J_1-q_1)|\nabla_\theta w|^2+(p J_2-q_2)|\Delta_\theta w|^2\\
&+(p J_3-q_3)|\nabla_\theta\Delta_\theta w|^2+(p-1)|\Delta^2_\theta w|^2\leq0.
\endaligned\ee
Notice that $p J_0-q_0=0$ is  equivalent to
\be\label{11gammaine}
p\frac{\Gamma(\frac{n}{2}-\frac{4}{p-1})\Gamma(4+\frac{4}{p-1})}{\Gamma(\frac{4}{p-1})\Gamma(\frac{n-8}{2}-\frac{4}{p-1})}
=\frac{\Gamma(\frac{n+8}{4})^2}{\Gamma(\frac{n-8}{4})^2},
\ee
see \cite{LWZ}. To solve this, let us the transform $k:=\frac{8}{p-1}$ and $k:=\frac{n-10}{2}-a$, then we can reduce the equation  $p J_0-q_0=0$   to the following
\be\nonumber\aligned
&n^4a^8+(-n^5-20n^3)a^6+(\frac{3}{8}n^6+5n^4+118n^2)a^4\\
&+(-\frac{1}{16}n^7+\frac{5}{4}n^5-2n^3-180n)a^2\\
&+81+\frac{1}{16}n^7-\frac{7}{16}n^6-2n^5+\frac{115}{8}n^4+16n^3-109n^2=0.
\endaligned\ee
We further let that $a^2=t$ then the above equation is reduced to a fourth-order equation.  Find  the roots, they are
Let
\be\label{LWZ-666} R_1(n)=\frac{n-10}{2}-d(n),\quad R_2(n)=\frac{n-10}{2}+d(n),\ee
where $d(n)$ is defined in \eqref{LWZ=333}.
Transform back
(note that $p>\frac{n+8}{n-8}$, hence $0<k<\frac{n-8}{2}$ and thus only root satisfy this),  we get the
so called Joseph-Lundgren exponent, see \eqref{12pcn}. We have the following lemma after a  direct calculation.


\bl\label{12dn1} For the $d(n)$ appears in the Joseph-Lundgren exponent, see \eqref{12pcn}, we have
\be\nonumber\aligned
&\lim_{n\rightarrow+\infty}\frac{d(n)}{\sqrt{n}}=1,\\
&d(n)<\sqrt{n}\;\;\hbox{for}\;\;n\geq18.
\endaligned\ee

\br
This Lemma gives optimal bound for $d(n)$ in the sense that if find the optimal constant $B$ satisfy
\be\nonumber\aligned
d(n)<B\sqrt{n} \;\;\hbox{for}\;\;n \geq n_0
\endaligned\ee
\er

\el


\bl\label{12pJ1} If $\frac{n-10}{2}-\sqrt{n}<k<\frac{n-10}{2}+\sqrt{n}$  for $n\geq118$, we have
$$p J_1-q_1>0.$$

\el
\bp Equivalently, we consider the function $W_1:=k(p J_1-q_1)$. From \eqref{12J1} and \eqref{12q}, we have that
\be\nonumber\aligned
 W_1&=-4k^7+(-128+12n)k^6+(-12n^2+324n-1696)k^5+(4n^3-264n^2\\
&+3488n-12032)k^4
+(68n^3-2168n^2+19056n-49216)k^3+(376n^3\\
&-8176n^2+55104n-115712)k^2
+(-144384-\frac{1}{16}n^6+\frac{3}{4}n^5+732n^3\\
&-13520n^2+78336n)k+384n^3-6912n^2+39936n-73728.
\endaligned\ee
Set $k=\frac{n-10}{2}+a\sqrt{n}, n=t^2$ (this is a key point). Hence,  by the assumption we have that $-1\leq a\leq1$. Thus
\be\nonumber\aligned
W_1&=-108+(\frac{1}{2}-\frac{3}{8}a^2)t^{12}+(\frac{3}{4}a^3-\frac{3}{4}a)t^{11}+(-\frac{51}{8}+\frac{3}{2}a^4+\frac{9}{4}a^2)t^{10}\\
&+(-3a^5+\frac{9}{4}a)t^9+(\frac{3}{4}-2a^6-9a^4+11a^2)t^8+(4a^7+2a^3+28a)t^7\\
&+(\frac{351}{2}+12a^6-2a^4-10a^2)t^6
+(-44a^5-32a^3-47a)t^5\\
&+(-132a^4-150a^2-157)t^4+(76a^3-224a)t^3+(228a^2-1126)t^2-36at.
\endaligned\ee
For the case of $0\leq a\leq1$, we get from the above identity that
\be\nonumber\aligned
W_1&\geq-108+\frac{1}{8}t^{12}-\frac{3}{4}t^{11}-\frac{51}{8}t^{10}-\frac{9}{4}t^9+\frac{3}{4}t^8
+\frac{351}{2}t^6-123t^5-439t^4\\
&\quad -224t^3-1126t^2-36t>0\;\;\hbox{if}\;\;n=t^2\geq117\;(t\geq10.7725).
\endaligned\ee
For the case of  $0\leq a\leq1$, we get  that
 \be\nonumber\aligned
W_1&\geq-108+\frac{1}{8}t^{12}-\frac{3}{4}t^{11}-\frac{51}{8}t^{10}-\frac{9}{4}t^9+\frac{3}{4}t^8
-34t^7+\frac{351}{2}t^6-439t^4\\
&-76t^3-1354t^2>0\;\;\hbox{if}\;\;n=t^2\geq118\;(t\geq10.8562627).
\endaligned\ee
\ep

\bl\label{12pJ2} If $\frac{n-10}{2}-\sqrt{n}<k<\frac{n-10}{2}+\sqrt{n}$  for $n\geq30$, we have
$$p J_2-q_2>0.$$

\el
\bp
Equivalently, we consider the function $W_2:=k(p J_2-q_2)$. From \eqref{12J2} and \eqref{12q}, we have that
 \be\nonumber\aligned
W_2=&6k^5+(144-12n)k^4+(6n^2-228n+1352)k^3+(84n^2-1544n\\
&+6272)k^2+(14464-\frac{3}{8}n^4+3n^3+338n^2-4512n)k+352n^2\\
&-4480n+13824.
\endaligned\ee
Set $k=\frac{n-10}{2}+a\sqrt{n}, n=t^2$,  hence  $-1\leq a\leq1$. It follows that
 \be\nonumber\aligned
W_2=&554+(3-\frac{3}{2}a^2)t^8+(3a^3-3a)t^7+(-\frac{51}{2}+3a^4+3a^2)t^6\\
&+(-6a^5-3a)t^5+(-6a^4+34a^2-51)t^4+(28a^3+80a)t^3\\
&+(92a^2+427)t^2+106at.
\endaligned\ee
For the case of $0\leq a\leq1$, from the above inequality we have that
\be\nonumber\aligned
W_2\geq&554+\frac{3}{2}t^8-3t^7-\frac{51}{2}t^6-9t^5-57t^4+427t^2\\
>&0\;\;\hbox{if}\;\;n=t^2\geq30\;(t\geq5.475795).
\endaligned\ee
For the case of  $-1\leq a\leq0$, we observe that
\be\nonumber\aligned
W_2\geq&554+\frac{3}{2}t^8-3t^7-\frac{51}{2}t^6-57t^4-108t^3+427t^2-106t\\
>&0\;\;\hbox{if}\;\;n=t^2\geq28\;(t\geq5.257108771).
\endaligned\ee

\ep


\bl\label{12pJ3} If $\frac{n-10}{2}-\sqrt{n}<k<\frac{n-10}{2}+\sqrt{n}$  for $n\geq11$, we have
$$p J_3-q_3>0.$$

\el

\bp
We consider the function $W_3:=k(p J_3-q_3)$. From \eqref{12J3} and \eqref{12q}, we have that
\be\nonumber\aligned
W_3=-4k^3+(-64+4n)k^2+(-n^2+48n-320)k-640+96n.
\endaligned\ee
Set $k=\frac{n-10}{2}+a\sqrt{n}, n=t^2$,  then $-1\leq a\leq1$.   Thus,
\be\nonumber\aligned
W_3=-140+(-2a^2+8)t^4+(4a^3-4a)t^3+(-4a^2-34)t^2-20at.
\endaligned\ee
Since $a\in[-1,1]$, from the above inequality we get that
\be\nonumber\aligned
W_3\geq-140+6t^4-4t^3-38t^2-20t>0\;\;\hbox{if}\;\;n=t^2\geq12\;(i.e., t\geq3.40511).
\endaligned\ee
\ep


For the lower dimension case, we can calculate numerically and  thus we have the following Lemma
\bl\label{12pJ4}
Consider the supercritical case $p>\frac{n+8}{n-8}$, i.e.,  $0<k<\frac{n-8}{2}$. We have the following facts.
\begin{itemize}
\item [(1)] If  $ 0<k<\frac{n-8}{2}$ and $n\leq 17$, then $pJ_1-q_1,pJ_2-q_2,pJ_3-q_3>0$;
\item [(2)] If  $ 18\leq n\leq 120$ and $R_1(n)<k<R_2(n)$, then $pJ_1-q_1>0$;
\item [(3)] If  $ 18\leq n\leq 29$ and $R_1(n)<k<R_2(n)$, then $pJ_2-q_2>0$;
\end{itemize}
where  $R_1(n), R_2(n)$ are given in \eqref{LWZ-666}.
\el

\bt\label{12stable} Let $u\in W^{4,2}_{loc}(\R\backslash \{0\})$ be a homogeneous and  stable solution of \eqref{12LE}
with $\frac{n+8}{n-8}<p<p_c(n)$. Assume that $|u|^{p+1}\in L^{p+1}(\R^n\{0\})$, then $u\equiv0$.

\et
\bp From the inequality \eqref{12JPQ}, combine with Lemmas \ref{12pJ1},\ref{12pJ2},\ref{12pJ3} and \ref{12pJ4}, we may get the
conclusion easily. \ep


\section{Energy estimates and Blow down analysis }

 In this section, we finish the  energy estimates for the solutions to \eqref{12LE}, which are important when we perform a
 blow-down analysis in the next section.

 \subsection{Energy estimates}

\bl\label{11lemma0} Let $u$ be a stable solution of \eqref{12LE}, then there exists a positive constant $C$ independent of $R$  such that

\be\label{11estimate}\aligned
&\int_{\R^n}|u|^{p+1}\eta^2+\int_{\R^{n}}|\Delta^2 u|^2\eta^2\\
&\leq C\big[\int_{\R^{n}}|\nabla\Delta u|^2|\frac{\nabla \eta^2}{\eta}|^2
+|\Delta u|^2\big((\frac{\Delta \eta^2}{\eta})^2\\
&\quad+(\frac{\nabla^2\eta^2}{\eta})^2\big)
+|\nabla u|^2(\frac{\nabla\Delta \eta^2}{\eta})^2+ u^2(\frac{\Delta^2 \eta^2}{\eta})^2\Big]
\endaligned\ee

\el
\bp
Multiply the equation \eqref{12LE} with $ u\eta^2$, where $\eta$ is a test function, we get that

\be\label{8eee}
\int_{\R^n}|u|^{p+1}\eta^2=\int_{\R^n}(\Delta^4u) u\eta^2=\int_{\R^n}(\Delta^2u)\Delta^2(u\eta^2).
\ee
Since $\Delta(\xi\eta)=\eta\Delta\xi+\xi\Delta\eta+2\nabla\xi\nabla\eta$, we have
\be\nonumber\aligned
\Delta^2 u\Delta^2(u\eta^2)&=(\Delta^2u)^2\eta^2+2\Delta^2u\Delta u\Delta\eta^2
+2\Delta^2 u\nabla\eta^2\nabla\Delta u\\
&\quad+u\Delta u \Delta^2\eta^2
+2\Delta^2u\nabla u\nabla\Delta\eta^2+2\Delta^2 u\Delta(\nabla u\nabla\eta^2)\\
&=(\Delta^2u)^2\eta^2+\hbox{the combination of terms with  lower order  than}(\Delta^2u)^2.
\endaligned\ee
Further
\be\label{10id}\aligned
|\Delta^2&(u\eta)|^2=(\Delta^2u)^2\eta^2+4(\Delta u)^2(\Delta\eta)^2+u^2(\Delta^2\eta)^2
+4(\nabla u\nabla\Delta\eta)^2\\
&\quad+4(\nabla\eta\nabla\Delta u)^2+4|\Delta(\nabla u\nabla\eta)|^2
+4\Delta^2u\eta\Delta u\Delta\eta+2(\Delta^2 u)\eta u\Delta^2\eta\\
&\quad+4(\Delta^2 u)\eta\nabla u\nabla\Delta\eta+4(\Delta^2u)\eta\nabla\eta\nabla\Delta u
+4(\Delta^2u) u\eta\Delta(\nabla u\nabla\eta)\\
&\quad+4\Delta u(\Delta\eta) u\Delta^2\eta
+8\Delta u\Delta\eta\nabla u\nabla\Delta\eta+8\Delta u\Delta\eta\nabla\eta\nabla\Delta u\\
&\quad+8\Delta u\Delta\eta\Delta(\nabla u\nabla\eta)
+4u\Delta^2\eta\nabla u\nabla\Delta\eta+4u\Delta^2\eta\nabla\eta\nabla\Delta u\\
&\quad+4u\Delta^2\eta\Delta(\nabla u\nabla\eta)+8(\nabla u\nabla\Delta\eta)(\nabla \eta\nabla\Delta u)\\
&\quad+8\nabla u\nabla\Delta\eta\Delta(\nabla u\nabla\eta)+8\nabla\eta\nabla\Delta u\Delta(\nabla u\nabla\eta)\\
&\quad=(\Delta^2u)^2\eta^2+\hbox{the combination of  terms with lower order than}(\Delta^2u)^2.
\endaligned\ee
On the other hand, by the stability condition, we have
\be\label{10stable1}
p\int_{\R^n}|u|^{p+1}\eta^2\leq\int_{\R^{n}}|\Delta^2(u\eta)|^2.
\ee
Combining with \eqref{8eee}, \eqref{10id} and \eqref{10stable1}, via the basic Cauchy inequality, we prove the lemma.
\ep


\vskip0.1in

\bl\label{12key}
Let $u$ be a stable solution of \eqref{12LE}. Then
\be\label{12key1}
\int_{B_R}|u|^{p+1}+\int_{B_R}|\Delta^2 u|^2\leq C R^{-8}\int_{B_{2R}}u^2,
\ee
\be\label{12key2}
\int_{B_R}|u|^{p+1}+\int_{B_R}|\Delta^2 u|^2\leq C R^{n-8\frac{p+1}{p-1}}.
\ee
\el

\bp
Now let $\eta=\xi^m$,  where $m>4$, $\xi=1$ in $B_{R/2}$ and $\xi=0$ in $B_R^cC$ satisfying $|\nabla \xi |\leq\frac{C}{R}$. Plug into the estimates in the previous lemma. Then
 \be\label{12F123}\aligned
&\int_{\R^n}|\Delta^2 u|^2\xi^{8m}+\int_{\R^n}|u|^{p+1}\xi^{8m}
\leq C\Big(\int_{\R^n}u^2g_0(\xi)+\int_{\R^n}|\nabla u|^2g_1(\xi)\\
&+\int_{\R^n}|\Delta u|^2g_2(\xi)
+\int_{\R^n}|\nabla\Delta u|^2g_3(\xi)\Big),
\endaligned\ee
where
\be\nonumber\aligned
&g_0(\xi):=\xi^{8m-8}\sum_{0\leq i+j+k+s+t+u+v+w=8}|\nabla^i\xi||\nabla^j\xi||\nabla^k\xi||\nabla^s\xi||\nabla^t\xi||\nabla^u\xi||\nabla^v\xi||\nabla^w\xi|,\\
&g_1(\xi):=\xi^{8m-6}\sum_{0\leq i+j+k+s+t+u=6}|\nabla^i\xi||\nabla^j\xi||\nabla^k\xi||\nabla^s\xi||\nabla^t\xi||\nabla^u\xi|,\\
&g_2(\xi):=\xi^{8m-4}\sum_{0\leq i+j+k+s=4}|\nabla^i\xi||\nabla^j\xi||\nabla^k\xi||\nabla^s\xi|,\\
&g_3(\xi):=\xi^{8m-2}\sum_{0\leq i+j=2}|\nabla^i\xi||\nabla^j\xi|,\\
\endaligned\ee
here $\nabla^0\xi:=\xi$ and notice that $g_m(\xi)\geq0$ for $m=0,1,2,3$.
Now we claim that
\be\nonumber\aligned
&g_1^2(\xi)\leq C g_0(\xi)g_2(\xi),\quad g_2^2(\xi)\leq C g_1(\xi)g_3(\xi),\quad g_3^2(\xi)\leq C \xi^{8m}g_2(\xi),\\
&|\nabla^2 g_3(\xi)|\leq C g_2(\xi), \quad |\nabla^2 g_2(\xi)|\leq C g_1(\xi), \quad|\nabla^2 g_1(\xi)|\leq C g_0(\xi).
\endaligned\ee
This claim can be checked directly and will be frequently applied to  our estimates below.
In the next, we evaluate  every term  in the right hand side of \eqref{12F123}.
By a   integrate by part, we have
\be\nonumber
|\nabla\Delta u|^2=\frac{1}{2}\Delta(\Delta u)^2-\Delta^2u\Delta u.
\ee
It follows that
\be\label{12d3}\aligned
&\int_{\R^n}|\nabla\Delta u|^2g_3(\xi)=\frac{1}{2}\int_{\R^n}\Delta(\Delta u)^2g_3(\xi)-\int_{\R^n}\Delta^2u\Delta u g_3(\xi)\\
&=\frac{1}{2}\int_{\R^n}(\Delta u)^2\Delta g_3(\xi)-\int_{\R^n}\Delta^2u\Delta u g_3(\xi)\\
&\leq \varepsilon_3\int_{\R^n}|\Delta^2 u|^2\xi^{8m}+C(\varepsilon_3)\int_{\R^n}|\Delta u|^2g_2(\xi),
\endaligned\ee
where $\varepsilon_3$ is a parameter to be determined later.
Integrating by part again  we have
\be\nonumber
(\Delta u)^2=\sum_{j,k}(u_j u_k)_{jk}-\sum_{j,k}(u_{jk})^2-2\nabla\Delta u\nabla u,
\ee
hence $(\Delta u)^2\leq\sum_{j,k}(u_j u_k)_{jk}-2\nabla\Delta u\nabla u.$ By which we get
\be\label{12d2}\aligned
\int_{\R^n}(\Delta u)^2g_2(\xi)&\leq \int_{\R^n}\sum_{j,k}(u_j u_k)_{jk}g_2(\xi)-2\int_{\R^n}\nabla\Delta u\nabla u g_2(\xi)\\
&=\int_{\R^n}\sum_{j,k}u_j u_kg_2(\xi)_{jk}-2\int_{\R^n}\nabla\Delta u\nabla u g_2(\xi)\\
&\leq C(\varepsilon_2)\int_{\R^n}|\nabla u|^2g_1(\xi)+\varepsilon_2\int_{\R^n}|\nabla\Delta u|^2g_3(\xi),
\endaligned\ee
where $\varepsilon_2$ is a parameter to be determined later.
From the differential identity, $|\nabla u|^2=\frac{1}{2}\Delta u^2-u\Delta u$, we get that
\be\label{12d1}\aligned
\int_{\R^n}|\nabla u|^2 g_1(\xi)&=\frac{1}{2}\int_{\R^n}\Delta u^2g_1(\xi)-\int_{\R^n}u\Delta u g_1(\xi)\\
&=\frac{1}{2}\int_{\R^n}u^2g_1(\xi)-\int_{\R^n}u\Delta u g_1(\xi)\\
&\leq C(\varepsilon_1)\int_{\R^n}u^2 g_0(\xi)+\varepsilon_1\int_{\R^n}(\Delta u)^2 g_2(\xi).
\endaligned\ee
Combining with \eqref{12d3}, \eqref{12d2} and \eqref{12d1}, by selecting the parameters  $\varepsilon_1,\varepsilon_2$
small enough, we can obtain that
\be\nonumber\aligned
&\int_{\R^n}|\nabla u|^2g_1(\xi)+\int_{\R^n}|\Delta u|^2g_2(\xi)+\int_{\R^n}|\nabla\Delta u|^2g_3(\xi)\\
&\leq C(\varepsilon_3)\int_{R^n}u^2g_0(\xi)+\varepsilon_3\int_{\R^n}|\Delta^2 u|^2\xi^{8m}.
\endaligned\ee
Combining the  above estimate     with \eqref{12F123} and   selecting $\varepsilon_3$ small enough, we have that
\be\nonumber\aligned
\int_{\R^n}|\Delta^2 u|^2\xi^{8m}+\int_{\R^n}|u|^{p+1}\xi^{8m}\leq C\int_{\R^n}u^2g_0(\xi).
\endaligned\ee
This proves \eqref{12key1}.
Further, we let $\xi=1$ in $B_R$ and $\xi=0$ in $B_{2R}^c$, satisfying $|\nabla\xi|\leq\frac{C}{R}$, then we have
\be\nonumber\aligned
&\int_{\R^n}|\Delta^2 u|^2\xi^{8m}+\int_{\R^n}|u|^{p+1}\xi^{8m}\leq C\int_{\R^n}u^2g_0(\xi)
\leq CR^{-8}\int_{\R^n}u^2\xi^{8m-8}\\
&\leq CR^{-8}\int_{\R^n}\big(|u|^{p+1}\xi^{(4m-4)(p+1)}\big)^{\frac{2}{p+1}}R^{n(1-\frac{2}{p+1})}.
\endaligned\ee
By choosing $m=\frac{p+1}{p-1}$, hence $(4m-4)(p+1)=8m$, it follows that \eqref{12key2} holds.

\ep


\subsection{Blow-down analysis and the proof of Theorem \ref{12th1}}

{\bf The proof of Theorem \ref{12th1}}. Firstly, we consider $1<p\frac{n+8}{n-8}$. If $p<\frac{n+8}{n-8}$, we can let
$R\rightarrow+\infty$ in \eqref{12key2} to get $u\equiv0$ directly.
  For $p=\frac{n+8}{n-8}$, hence $n=8\frac{p+1}{p-1}$, \eqref{12key2} gives that
\be\nonumber
\int_{\R^n}|\Delta^2 u|^2+|u|^{p+1}<+\infty.
\ee
Hence
\be\nonumber
\lim_{R\rightarrow+\infty}\int_{B_{2R}(x)\setminus B_R(x)}|\Delta^2 u|^2+|u|^{p+1}=0.
\ee
Then by Lemma \ref{12key} and noting that now $n=8\frac{p+1}{p-1}$, we have
\be\nonumber\aligned
&\int_{B_R(x)}|\Delta^2 u|^2+|u|^{p+1}\leq CR^{-8}\int_{B_{2R}\setminus B_R(x)}u^2\\
&\leq CR^{-8}\Big(\int_{B_{2R}\setminus B_R(x)}|u|^{p+1}\Big)^{\frac{2}{p+1}}R^{n(1-\frac{2}{p+1})}
\leq\Big(\int_{B_{2R}\setminus B_R(x)}|u|^{p+1}\Big)^{\frac{2}{p+1}},
\endaligned\ee
letting $R\rightarrow+\infty$, we get that $u\equiv0$.
\vskip0.1in

Secondly, we consider the supercritical case, i.e., $p>\frac{n+8}{n-8}$. We divide  the proof into several steps.

\vskip0.2in
{\bf Step 1: $\lim_{\la \rightarrow+\infty}E(u,0,\la)<+\infty.$}
From Theorem \ref{12mono2}, we know that $E$ is nondecreasing w.r.t. $\la$, note that

\be\nonumber\aligned
E(u,0,\la)\leq\frac{1}{\la}\int_{\la}^{2\la}E(u,0,t)dt\leq\frac{1}{\la^2}\int_{\la}^{2\la}\int_{t}^{t+\la}E(u,0,\gamma)d\gamma dt,
\endaligned\ee
where $C>0$ is constant independent of $\gamma$.
From Lemma \ref{12key}, we have that
\be\nonumber\aligned
\frac{1}{\la^2}&\int_\la^{2\la}\int_{t}^{t+\la}\gamma^{8\frac{p+1}{p-1}-n}\big[\int_{B_\gamma}\frac{1}{2}|\Delta^2 u|^2dx -\frac{1}{p+1}\int_{B_\gamma}|u|^{p+1}dx\big]d\gamma dt\leq C,
\endaligned\ee
where $C>0$ is independent of $\gamma$. Further,
\be\label{12type1}\aligned
&\frac{1}{\la^2} \int_\la^{2\la}\int_{t}^{t+\la}\int_{\pa B_\gamma}\gamma^{8\frac{p+1}{p-1}-n-7}\Big[C_0u
+C_1\gamma\pa_r u
+C_2\gamma^2 \pa_{rr} u+C_3\gamma^3\pa_{rrr} u\Big]\\
&\quad\quad\Big[C_0^1u+C_1^1\gamma\pa_r u+C_2^1\gamma^2\pa_{rr}u\Big]\\
&\leq C\frac{1}{\la^2}\int_\la^{2\la}\int_{t}^{t+\la}t^{8\frac{p+1}{p-1}-n-8}\int_{\pa B_\gamma} \big[u^2+\gamma^2(\pa_r u)^2+\gamma^4(\pa_{rr}u)^2+\gamma^6(\pa_{rrr}u)^2\big]\\
&\leq C\frac{1}{\la^2}\int_\la^{2\la}t^{8\frac{p+1}{p-1}-n-8}\int_{B_{3\la}}\big[u^2+\gamma^2(\pa_r u)^2+\gamma^4(\pa_{rr}u)^2+\gamma^6(\pa_{rrr}u)^2\big]\\
&\leq C \la^{n-8\frac{p+1}{p-1}+8}\frac{1}{\la^2}\int_\la^{2\la}t^{8\frac{p+1}{p-1}-n-7}dt\\
&\leq C.
\endaligned\ee
Integrating by part if necessary, the remaining    terms can be treated similarly as the estimate of \eqref{12type1}.

\vskip0.1in
{\bf Step 2:  } For any $\la>0$, recall  the definition

\be\nonumber
u^\la(x):=\la^{\frac{8}{p-1}}u(\la x)
\ee
and $u^\la$ is also a smooth stable solution of \eqref{12LE} in $\R^n$.
By rescaling the estimate \eqref{12key2} in Lemma \ref{12key}, for any $\la>0$ and balls $B_r(x)\subset \R^n$, we have that
\be\nonumber
\int_{B_r(x)} |\Delta^2 u^\la|^2+|u^\la|^{p+1}\leq Cr^{n-8\frac{p+1}{p-1}}.
\ee
In particular, $u^\la$ are uniformly bounded in $L^{p+1}_{loc}(\R^n)$ and $\Delta^2 u^\la$ are uniformly bounded in $L^2_{loc}(\R^n)$.
By elliptic estimates, $u^\la$ are also uniformly bounded in $W^{4,2}_{loc}(\R^n)$. Hence, up to s sequence of $\la\rightarrow+\infty$,
we can assume that $u^\la\rightarrow u^\infty$ weakly in $W^{4,2}_{loc}\cap L^{p+1}_{loc}(\R^n)$. By the Sobolev embedding,
$u^\la\rightarrow u^\infty$ in $W^{3,2}_{loc}(\R^n)$. Then for any ball $B_R(0)$, by the interpolation theorem and   the estimate \eqref{12key2},
for any $q\in[1,p+1)$ as $\la\rightarrow+\infty$, we obtain that
\be\label{12q12}
\parallel u^\la-u^\infty\parallel_{L^q(B_R(0))}\leq\parallel u^\la-u^\infty\parallel_{L^1(B_R(0))}^t\parallel u^\la-u^\infty\parallel_{L^{p+1}(B_R(0))}^{1-t},
\ee
where $t\in(0,1]$ satisfying  $\frac{1}{q}=t+\frac{1-t}{q+1}$.
That is,  $u^\la\rightarrow u^\infty$ in $L^{q+1}_{loc}(\R^n)$ for any $q\in [1,p+1)$.
For any $\varphi\in C_0^\infty(\R^n)$, we have that
\be\nonumber\aligned
&\int_{\R^n}\Delta^2 u^\infty \Delta^2 \varphi-|u^\infty|^{p-1}u^\infty\varphi
=\lim_{\la\rightarrow+\infty}\int_{\R^n}\Delta^2 u^\la\Delta^2 \varphi-|u^\la|^{p-1}u^\la\varphi,\\
&\int_{\R^n}|\Delta^2 \varphi|^2-p|u^\infty|^{p-1}\varphi^2=\lim_{\la\rightarrow+\infty}\int_{\R^n}
|\Delta^2 \varphi|^2-p|u^\la|^{p-1}\varphi^2.
\endaligned\ee
Therefore,  $u^\infty\in W^{4,2}_{loc}\cap L^{p+1}_{loc}(\R^n)$ is  a stable solution of \eqref{12LE} in $\R^n$.

\vskip0.1in

{\bf Step 3: } We claim that the function $u^\infty$ is homogeneous.
Due to the scaling invariance of the functional $E$ (i.e., $E(u,0,R\la)=E(u^\la,0,R)$) and the monotonicity formula, for any given
$R_2>R_1>0$, we have that
\be\nonumber\aligned
0&=\lim_{i\rightarrow+\infty}\Big(E(u,0,R_2\la_i)-E(u,0,R_1\la_i)\Big)\\
&=\lim_{i\rightarrow+\infty}\Big(E(u^{\la_i},0,R_2)-E(u^{\la_i},0,R_1)\Big)\\
&\geq C(n,p)\liminf_{i\rightarrow+\infty}\int_{B_{R_2}\setminus B_{R_1}}r^{8\frac{p+1}{p-1}-n-8}\big(\frac{8}{p-1}u^{\la_i}+r\frac{\pa u^{\la_i}}{\pa r}\big)^2\\
&\geq C(n,p)\int_{B_{R_2}\setminus B_{R_1}}r^{8\frac{p+1}{p-1}-n-8}\big(\frac{8}{p-1}u^{\infty}+r\frac{\pa u^{\infty}}{\pa r}\big)^2.\\
\endaligned\ee
In the last inequality we have used the weak convergence of the sequence $(u^{\la_i})$ to the function $u^\infty$
in $W^{1,2}_{loc}(\R^n)$ as $i\rightarrow+\infty$.
This equality above implies that
\be\nonumber\aligned
\frac{8}{p-1}\frac{u^\infty}{r}+\frac{\pa u^\infty}{\pa r}=0,\;\;\hbox{a.e.}\;\;\hbox{in}\;\;\R^n.
\endaligned\ee
Integrating in $r$ shows that
\be\nonumber
u^\infty(x)=|x|^{-\frac{8}{p-1}} u^\infty(\frac{x}{|x|}).
\ee
That is, $u^\infty$ is homogeneous.
\vskip0.1in
{\bf Step 4: $u^\infty=0$. } This is a direct consequence of Theorem \ref{12stable} since $u^\infty$ is homogeneous.
Since this holds for the limit of any sequence $\la\rightarrow+\infty$, by \eqref{12q12} we get that
\be\nonumber
\lim_{\la\rightarrow+\infty} u^\la\;\;\hbox{strongly in}\;\;L^2(B_4(0)).
\ee
\vskip0.1in
{\bf Step 5: $u\equiv0$}. For all $\la\rightarrow+\infty$, we see that
\be\nonumber
\lim_{\la\rightarrow+\infty} \int_{B_4(0)}(u^\la)^2=0.
\ee
By \eqref{12key1} in Lemma \ref{12key}, we have that
\be\label{12e517}
\lim_{\la\rightarrow+\infty}\int_{B_3(0)}|\Delta^2 u^\la|^2+|u^\la|^{p+1}\leq\lim_{\la\rightarrow+\infty}\int_{B_4(0)}(u^\la)^2=0.
\ee
By the elliptic interior $L^2-$ estimate, we get that
\be\nonumber
\lim_{\la\rightarrow+\infty}\int_{B_2(0)}\sum_{j\leq4}|\nabla^j u^\la|^2=0.
\ee
In particular, we can choose a sequence $\la_i\rightarrow+\infty$ such that
\be\nonumber
\int_{B_2(0)}\sum_{j\leq4}|\nabla^j u^{\la_i}|^2\leq 2^{-i}.
\ee
Hence we have
\be\nonumber\aligned
\int_{1}^2\sum_{i=1}^{+\infty}\int_{\pa B_r}\sum_{j\leq4}|\nabla^j u^{\la_i}|^2dr
\leq\sum_{i=1}^{+\infty}\int_1^2\int_{\pa B_r}\sum_{j\leq4}|\nabla^j u^{\la_i}|^2\leq1.
\endaligned\ee
Therefore, the function
\be\nonumber\aligned
g(r):=\sum_{i=1}^\infty\int_{\pa B_r}\sum_{j\leq4}|\nabla^j u^{\la_i}|^2\in L^1(1,2).
\endaligned\ee
Then there exists an $r_0\in (1,2)$ such that $g_{(r_0)}<+\infty$, by which we get that
\be\nonumber\aligned
\lim_{i\rightarrow+\infty}\|u^{\la_i}\|_{W^{4,2}(\pa B_{r_0})}=0.
\endaligned\ee
Combine with \eqref{12e517} and the scaling invariance of $E(u,0,\la)$, we have
\be\nonumber\aligned
\lim_{i\rightarrow+\infty}E(\la r_0,0,u)=\lim_{i\rightarrow+\infty}E( r_0,0,u^{\la_i})=0.
\endaligned\ee
Since $\la_i r_0\rightarrow +\infty$ and $E(r,0,u)$ is nondecreasing in $r$, we get that
\be\nonumber\aligned
\lim_{i\rightarrow+\infty}E(\la r_0,0,u)=0.
\endaligned\ee
By the smoothness of $u$, $\lim_{i\rightarrow0}E(\la r_0,0,u)=0.$ Again by the monotonicity
of $E(r,0,u)$ and Step 4, we obtain that
\be\nonumber\aligned
E(r,0,u)=0\;\;\hbox{for all}\;\;r>0.
\endaligned\ee
Therefore by the monotonicity formula (i.e., Theorem \ref{12mono2}) we known that $u$ is homogeneous, then
$u\equiv0$ by Theorem \ref{12stable}.   \hfill$\Box$


\section {Finite Morse index solution}

In this section, we prove Theorem \ref{12th2}. Firstly, we have
\bl\label{12lemma1}
Let $u$ be a smooth (positive or sing changing) solution of \eqref{12LE} with finite Morse index, then there exist  constant $C>0$
and $R_0>0$ such that
\be\nonumber
|u(x)|\leq C|x|^{-\frac{8}{p-1}},\;\;\hbox{for any}\;\;x\in B_{R_0}^c.
\ee

\el

\bp Assume that $u$ is stable outside $B_{R_0}^c$. For any $x\in B_{R_0}^c$, let $M(x):=|u(x)|^{\frac{p-1}{8}}$
and $d(x)=|x|-R_0$.
Assume that the conclusion does  not holds, then there exists a sequence of $x_k\in B_{R_0}^c$ such that
\be\nonumber
M(x_k)d(x_k)\geq 2k.
\ee
Since $u$ is bounded on any compact set of $\R^n$, $d(x_k)\rightarrow +\infty$.
By the doubling Lemma (see \cite{Souplet2007}), there exists  another sequence $y_k\in B_{R_0}^c$ such that
\be\nonumber\aligned
&M(y_k)d(y_k)\geq2k,\;M(y_k)\geq M(x_k),\\
&M(z)\leq 2M(y_k)\;\;\hbox{for any}\;\;z\in B_{R_0}^c\;\;\hbox{such that }\;\;|z-y_k|\leq\frac{k}{M(y_k)}.
\endaligned\ee
Now we define
\be\nonumber
u_k(x):=M(y_k)^{-\frac{8}{p-1}}u(y_k+M(y_k)^{-1}x)\;\;\hbox{for}\;\;x\in B_R(0).
\ee
This and the above arguments give that
$u_k(0)=1,|u_k|\leq2^{\frac{8}{p-1}}$in $B_R(0)$. Further, $B_{k/M(y_k)}\cap B_{R_0}=\emptyset$, which implies that
$u$ is a stable solution in $B_{k/M(y_k)}(y_k)$, hence $u_k$ is stable in $B_R(0)$.
By elliptic regularity theory, $u_k$ are uniformly bounded in $C^9_{loc}(\R^n)$, up to s subsequence,
$u_k$ convergent to $u_\infty$ in $C^8_{loc}(\R^n)$. By the above conclusions on $u_k$, we have
\begin{itemize}
\item [(1)]   $|u_\infty(0)|=1$;
\item [(2)]  $|u_\infty|\leq 2^{\frac{8}{p-1}}$ in $\R^n$;
\item[(3)]  $u_\infty$ is a smooth stable solution of \eqref{12LE} in $\R^n$.
\end{itemize}
By the Liouville theorem for stable solution, i.e., Theorem \ref{12th1}, we get that
$u_\infty\equiv0$, this is a contradiction.
\ep


\bc\label{12coro1}
Under the same assumptions in the above Lemma \ref{12lemma1}, there exist constant $C>0$ and $R_0$ such that for all $x\in B_{R_0}^c$,
\be\nonumber
\sum_{0\leq j\leq7}|x|^{\frac{8}{p-1}+j}|\nabla^j u(x)|\leq C.
\ee
\ec
\bp
For any $x_0$ with $|x_0|>R_0$, take $\la=\frac{|x_0|}{2}$ and define
\be\nonumber
\overline{u}(x):=\la^{\frac{8}{p-1}}u(x_0+\la x).
\ee
By the previous Lemma, $\overline{u}(x)\leq C$ in $B_1(0)$.
By the elliptic regularity theory we have
\be\nonumber
\sum_{0\leq j\leq7}|\nabla^k \overline{u}(0)|\leq C.
\ee
Scaling back we get the conclusion. \ep


\subsection{The proof of Theorem \ref{12th2}-(1)}

This   is about the subcritical case, i.e.,  $1<p<\frac{n+8}{n-8}$. Firstly, we cite the following Pohozaev identity (see \cite{Wei1999}).

\bl For any function $u$ satisfying \eqref{12LE}, there holds
\be\nonumber
(\frac{n-8}{2}-\frac{n}{p+1})\int_{B_R}|u|^{p+1}=\int_{\pa B_R}B_4(u)d\sigma,
\ee
where
\be\nonumber\aligned
&B_4(u)=(2-\frac{n}{2})\sum_{k=1}^2(-\Delta)^{4-k}u\frac{\pa (-\Delta u)^{k-1}}{\pa n}-\frac{R}{p+1}|u|^{p+1}\\
&-(2-\frac{n}{2})\sum_{k=1}^2\frac{\pa (-\Delta )^{4-k}u}{\pa n}(-\Delta )^{k-1}u+\frac{1}{2}(-\Delta )^4uR
+2(-\Delta )^3u\frac{\pa u}{\pa n}\\
&-2\frac{\pa (-\Delta)^3u}{\pa n}u
+\sum_{k=1}^2<x,\nabla(-\Delta )^{k-1}u>\frac{\pa (-\Delta)^{4-k}u}{\pa n}\\
&-\sum_{k=1}^2(-\Delta)^{4-k}u\frac{\pa <x,\nabla(-\Delta)^ku>}{\pa n}.
\endaligned\ee
\el

\noindent {\bf The proof of Theorem \ref{12th2}-(1)}.
By Corollary \ref{12coro1}, for any $R>R_0$, noting that $p<\frac{n+8}{n-8}$ (hence $n-8\frac{p+1}{p-1}<0$), we have the following estimate:
\be\nonumber\aligned
\int_{\pa B_R}|B_4(u)|d\sigma\leq C\int_{\pa B_R}R^{-\frac{16}{p-1}-7}d\sigma
\leq C R^{n-8\frac{p+1}{p-1}}\rightarrow0\;\;\hbox{as}\;\;R\rightarrow+\infty.
\endaligned\ee
Letting $R\rightarrow+\infty$ in the above Pohozaev identity, we get that
\be\nonumber\aligned
(\frac{n-8}{2}-\frac{n}{p+1})\int_{\R^n}|u|^{p+1}=0.
\endaligned\ee
Since $\frac{n-8}{2}-\frac{n}{p+1}<0$, we see that $u\equiv0$.  \hfill $\Box$


\subsection{The proof of Theorem \ref{12th2}-(3)}

Recall the assumption   $p=\frac{n+8}{n-8}$ (critical case) in  Theorem \ref{12th2}-(3).
Since $u$ is stable outside $B_{R_0}$, Lemma \ref{12key} holds if the support of $\eta$ is outside $B_{R_0}$.
Take $\varphi\in C_0^\infty(B_{2R_0}\setminus B_{R_0})$ such that $\varphi=1$ in $B_{R}\setminus B_{3R_0}$
and $\sum_{0\leq j\leq7}|x|^j|\nabla^j u|\leq1000$.
Then by choosing $\eta=\varphi^m$, where $m$ is bigger than 1, we get that
\be\nonumber\aligned
\int_{B_R\setminus B_{3R_0}}|\Delta^2 u|^2+|u|^{p+1}\leq C.
\endaligned\ee
Letting $R\rightarrow+\infty$, we have
\be\label{12finite}\aligned
\int_{\R^n}|\Delta^2 u|^2+|u|^{p+1}<+\infty.
\endaligned\ee
By the interior elliptic estimates and  the Holder's inequality, we have
\be\nonumber\aligned
&R^{-6}\int_{B_{2R}\setminus B_R}|\nabla u|^2\leq C\int_{B_{3R}\setminus B_{R/2}}|\Delta^2 u|^2+C\Big(\int_{B_{3R}\setminus B_{R/2}}|u|^{p+1}\Big)^{\frac{2}{p+1}},\\
&R^{-4}\int_{B_{2R}\setminus B_R}|\Delta u|^2\leq C\int_{B_{3R}\setminus B_{R/2}}|\Delta^2 u|^2+C\Big(\int_{B_{3R}\setminus B_{R/2}}|u|^{p+1}\Big)^{\frac{2}{p+1}},\\
&R^{-2}\int_{B_{2R}\setminus B_R}|\nabla\Delta u|^2\leq C\int_{B_{3R}\setminus B_{R/2}}|\Delta^2 u|^2+C\Big(\int_{B_{3R}\setminus B_{R/2}}|u|^{p+1}\Big)^{\frac{2}{p+1}},\\
&R^{-8}\int_{B_{2R}\setminus B_R}| u|^2\leq C\int_{B_{3R}\setminus B_{R/2}}|\Delta^2 u|^2+C\Big(\int_{B_{3R}\setminus B_{R/2}}|u|^{p+1}\Big)^{\frac{2}{p+1}},\\
\endaligned\ee
where $C$ is a universal constant independent of $R$.
Therefore, we have that
\be\nonumber\aligned
&\max\Big\{R^{-6}\int_{B_{2R}\setminus B_R}|\nabla u|^2,\;\; R^{-4}\int_{B_{2R}\setminus B_R}|\Delta u|^2,\;\; R^{-2}\int_{B_{2R}\setminus B_R}|\nabla\Delta u|^2,\\
&R^{-8}\int_{B_{2R}\setminus B_R}| u|^2\Big\}\rightarrow0\;\;\hbox{as}\;\;R\rightarrow+\infty.
\endaligned\ee
On the other hand, testing \eqref{12LE} with $u\eta^2$, we get that
\be\nonumber\aligned
\int_{\R^n}|\Delta^2 u|^2\eta^2-|u|^{p+1}\eta^2=-\int_{\R^n}\Delta^2 u\Big(\Delta^2(u\eta^2)-\Delta^2u\eta^2\Big)
\endaligned\ee
and
\be\nonumber\aligned
\Delta^2(u\eta^2)-\Delta^2u\eta^2&=2\Delta u\Delta\eta^2+2\nabla\Delta u\nabla\eta^2\\
&+u\Delta^2\eta^2+2\nabla u\nabla\Delta \eta^2+2\Delta(\nabla u\nabla\eta^2).
\endaligned\ee
Notice that the highest order derivative about $u$ of the  above expression if $\nabla\Delta u$.
By selecting $\eta(x)=\xi(\frac{x}{R})^{4m},m>1$ and $\xi\in C_0^\infty(B_2)$, $\xi=1$ in $B_1$
and $\sum_{1\leq j\leq5}|\nabla^j u|\leq1000$, we get that
\be\nonumber\aligned
&\Big|\int_{\R^n}|\Delta^2 u|^2\xi(\frac{x}{R})^{8m}-|u|^{p+1}\xi(\frac{x}{R})^{8m}\Big|
\leq C\Big(
R^{-6}\int_{B_{2R}\setminus B_R}|\nabla u|^2\\
&+R^{-4}\int_{B_{2R}\setminus B_R}|\Delta u|^2
+R^{-2}\int_{B_{2R}\setminus B_R}|\nabla\Delta u|^2+R^{-8}\int_{B_{2R}\setminus B_R}| u|^2
\Big).
\endaligned\ee
Now letting $R\rightarrow+\infty$, we obtain that
\be\nonumber\aligned
\int_{\R^n}|\Delta^2 u|^2-|u|^{p+1}=0,
\endaligned\ee
Combining this  with \eqref{12finite} we get the conclusion.\hfill $\Box$


\subsection{The proof of Theorem \ref{12th2}-(2)}

This is  the   supercritical case:    $p>\frac{n+8}{n-8}$.   Firstly, we have
\bl
There exists a constant $C>0$ such that $E(r,0,u)\leq C$ for all $r>R_0$.
\el
\bp
From the monotonicity formula, combine with the derivative estimate, i.e., Corollary \ref{12coro1}, we have the following
estimates:
\be\nonumber\aligned
&E(r,0,u)\leq Cr^{8\frac{p+1}{p-1}-n}\Big(|\Delta^2 u|^2+|u|^{p+1}\Big)\\
&+C\Big(\sum_{s,t\leq5,s+t\leq7}r^{8\frac{p+1}{p-1}-n-7+s+t}\int_{\pa B_r}|\nabla^s u||\nabla^t u|\\
&\leq C.
\endaligned\ee
\ep
As a consequence, we have the following
\bc
\be\nonumber
\int_{B_{3R_0}^c}\frac{\Big(\frac{8}{p-1}u(x)+|x|\frac{\pa u(x)}{\pa r}\Big)^2}{|x|^{n-8\frac{p+1}{p-1}}}dx<+\infty.\ee
\ec
As before, we define the blowing down sequence,
\be\nonumber\aligned
u^\la(x):=\la^{\frac{8}{p-1}}u(\la x).
\endaligned\ee
By Lemma \ref{12lemma1} and Corollary \ref{12coro1},  we know  that $u^\la$ are uniformly bounded in \\ $C^9(B_r(0)\backslash B_{1/r}(0))$
for any fixed $r>1$ and moreover, $u^\la$ is stable outside $B_{R_0/\la}$. And there exists a function $u^\infty\in C^8(\R^n\backslash\{0\})$,
such that up to a subsequence of $\la\rightarrow+\infty$, $u^\la$ convergent to $u^\infty$
in  $C^8(\R^n\backslash\{0\})$, $u^\infty$  is a stable solution of \eqref{12LE} in $\R^n\{0\}$. For any $r>1$, by the previous Corollary, we have
\be\nonumber\aligned
&\int_{B_r\setminus B_{1/r}}\frac{\Big(\frac{8}{p-1}u^\infty(x)+|x|\frac{\pa u^\infty(x)}{\pa r}\Big)^2}{|x|^{n-8\frac{p+1}{p-1}}}dx\\
&=\lim_{\la\rightarrow+\infty}\int_{B_r\setminus B_{1/r}}\frac{\Big(\frac{8}{p-1}u^\la(x)+|x|\frac{\pa u^\la(x)}{\pa r}\Big)^2}{|x|^{n-8\frac{p+1}{p-1}}}dx\\
&=\lim_{\la\rightarrow+\infty}\int_{B_{\la r}\setminus B_{\la/r}}\frac{\Big(\frac{8}{p-1}u(x)+|x|\frac{\pa u(x)}{\pa r}\Big)^2}{|x|^{n-8\frac{p+1}{p-1}}}dx\\
&=0.
\endaligned\ee
Hence,
\be\nonumber\aligned
\frac{8}{p-1}u^\infty(x)+|x|\frac{\pa u^\infty(x)}{\pa r}=0\;\;a.e.,
\endaligned\ee
that is,  $u$ is homogeneous, by Theorem \ref{12stable}, we get that $u^\infty\equiv0$ if $p<p_c(n)$.
Since this holds for any limit of $u^\la$ as $\la\rightarrow+\infty$, then we have
\be\nonumber\aligned
\lim_{|x|\rightarrow+\infty}|x|^{\frac{8}{p-1}}|u(x)|=0.
\endaligned\ee
Then as the proof of Corollary \ref{12coro1}, we have
\be\nonumber\aligned
\lim_{|x|\rightarrow+\infty}\sum_{0\leq j\leq7}|x|^{\frac{8}{p-1}+j}|\nabla^j u(x)|=0.
\endaligned\ee
Therefor, for any $\varepsilon>0$, take an $R_0$ such that for $|x|>R_0$, there holds
\be\nonumber\aligned
\sum_{0\leq j\leq7}|x|^{\frac{8}{p-1}+j}|\nabla^j u(x)|\leq \varepsilon.
\endaligned\ee
Then for any $r\gg R_0$, we have
\be\nonumber\aligned
E(r,0,u)& \leq  C r^{8\frac{p+1}{p-1}-n}\int_{B_R(0)}|\Delta^2 u|^2+|u|^{p+1}\\
&\quad +C\varepsilon r^{8\frac{p+1}{p-1}-n}\Big(\int_{B_r(0)\setminus B_{R_0}(0)}|x|^{-8\frac{p+1}{p-1}}+r\int_{\pa B_r(0)}|x|^{-8\frac{p+1}{p-1}}\Big)\\
&\leq C(R_0)\Big(r^{8\frac{p+1}{p-1}-n}+\varepsilon\Big).
\endaligned\ee
Therefore,  we obtain that
\be\nonumber\aligned
\lim_{r\rightarrow+\infty} E(r,0,u)=0
\endaligned\ee
since $8\frac{p+1}{p-1}+1-n<0$ and $\varepsilon$ can be arbitrarily small.
On the other hand, since $u$ is smooth, we have $\lim_{r\rightarrow0} E(r,0,u)=0$. Thus
$E(r,0,u)=0$ for all $r>0$. By the monotonicity formula we get that $u$ is homogeneous and hence by
Theorem \ref{12stable}, we derive that $u\equiv0$.  This completes the proof. \hfill $\Box$




\begin{thebibliography}{10}

\bibitem{Caffarelli1989}
L. Caffarelli, B. Gidas, J. Spruck, Asymptotic symmetry and local behavior of semilinear elliptic equations with
critical sobolev growth, {\it Comm. Pure Appl. Math.,} {\bf 42}, no. 3, 271-297, 1989.

\bibitem{Chang2011}
Sun-Yung Alice Chang and Maria del Mar Gonzalez, Fractional laplacian in conformal geometry,
{\it  Advances in Mathematic, } {\bf  226}, no. 2, 1410-1432, 2011.







\bibitem{Gidas-Sp=1981}  B. Gidas, J. Spruck, Global and local behavior of positive solutions of nonlinear elliptic equations, {\it Comm. Pure
Appl. Math.}{\bf  34}(1981), 525-598.



\bibitem{Gidas-Sp=1981-2} B. Gidas, J. Spruck, A priori bounds for positive solutions of nonlinear elliptic equations, {\it Comm. Partial Differential
Equations}, {\bf  6}(1981), 883-901.



\bibitem{Wei0=1}
J. Davila, L. Dupaigne, J. Wei, on the fractional Lane-Emden equation, Trans. Amer. Math. Soc. accepted for publication.

\bibitem{Wei=2}
J. Davila, L. Dupaigne, and K. Wang, J. Wei, A monotonicity formula and a liouville-type Theorem
 for a fourth order supercritical problem, {\it Advances in Mathematicas,}  {\bf 258}, 240-285, 2014.






\bibitem{Farina2005} A. Farina,   Liouville-type results for solutions of $-\Delta u=|u|^{p-1}u$ on unbounded domains of $\R^N,$
 {\it C. R. Math. Acad. Sci. Paris, }{\bf  341}, no. 7, 415-418, 2005




\bibitem{Farina2007}
A. Farina, On the classification of soultions of the Lane-Emden equation on unbounded domains of $\R^N$,
{\it J. Math. Pures Appl., } {\bf 87}(9),    no. 5, 537-561, 2007.




\bibitem{Fowler=1} R. Fowler, The solution of Emden's and similar differential equations, {\it Monthly Notices Roy.
Astronom. Soc., }{\bf 91}(1930), 63-91.

\bibitem{Fowler=2} R. Fowler,  Further studies of Emden's and similar differential equations, {\it Quarterly J. Math.
Oxford Ser.} {\bf  2}(1931), 259-288.








\bibitem{Gazzola2006}
F. Gazzola, H. C. Grunau, Radial entire solutions of supercritical biharmonic equations, {\it Math. Annal., } {\bf 334}, 905-936, 2006.








\bibitem{Wei1=2}
M. Fazly, J. Wei, on finite morse index solutions of higher order fractional Lane-Emden equations, American Journal of Mathematics accepted for publication.



\bibitem{Joseph1972}
D. D. Joseph, T. S. Lundgren, Quasilinear Dirichlet problems driven by positive sources,
{\it Arch. Rational Mech. Anal., } {\bf 49}, 241-269, 1972/73.




\bibitem{Lin1998}
C. S. Lin, A classification of soluitions of a conformally invariant fourth order equation in $\R^N$,
{\it Comment. Math. Helv., }  {\bf 73}, 206-231, 1998.

\bibitem{Mitidieri-96} E. Mitidieri, Nonexistence of positive solutions of semilinear elliptic systems in $\R^N$,
{\it Differential Integral Equations,}{\bf 9}(1996),  465-479.





\bibitem{LWZ2016=3} S. Luo, J. Wei and W. Zou, On the  Triharmonic  Lane-Emden Equation, submitted 2016.

\bibitem{LWZ} S. Luo, J. Wei and W. Zou,  On a transcendental equation involving quotients of Gamma functions,
 Proc. Amer. Math. Soc. accepted for publication.




\bibitem{Souplet2007}
P. Polacik, P. Quittner, P. Souplet, Singularity and decay estimates in superlinear problems via Liouville-type theorems, I. Elliptic equations and systems, {\it Duke Math. J.} {\bf 139}, no. 3, 555-579,  2007


\bibitem{QS2007} P. Quittner and Ph. Souplet, {\it Superlinear Parabolic Problems: Blow-up, Global Existence and Steady States}, Birkhausser, 2007, ISBN 978-3-7643-8442-5.




\bibitem{Serrin-Zou-1998} J. Serrin, H. Zou,  Existence of positive solutions of the Lane-Emden system. Dedicated to Prof. C. Vinti (Italian) (Perugia, 1996). Atti Sem. Mat. Fis. Univ. Modena, {\bf 46}(1998), suppl., 369-380.



\bibitem{Souplet2009} P. Souplet,  The proof of the Lane-Emden conjecture in four space dimensions, {\it  Adv. Math.,}  221 (2009), no. 5, 1409-1427.

\bibitem{WangXF-98-TAMS}X. F. Wang,  On the Cauchy problem for reaction-diffusion equations, {\it  Trans. Amer. Math. Soc.} {\bf  337}(1993), no. 2, 549-590.


\bibitem{Wei1999}
J. Wei, X. Xu, Classification of solutions of higher order conformally invariant equations,
{\it Math. Ann., } {\bf  313}, no. 2, 207-228,  1999.




\end{thebibliography}

\end{document}